\documentclass[12pt,leqno]{amsart}

\usepackage{amscd}
\usepackage{amsmath}
\usepackage{amsfonts}
\usepackage{amssymb}
\usepackage{url}
\usepackage{graphicx}
\usepackage{pictexwd,dcpic}
\usepackage{dsfont}

\title{On realizability of $p$-groups as Galois groups}
\author{Ivo M. Michailov, Nikola P. Ziapkov}
\address{Faculty of Mathematics and Informatics, Shumen University "Episkop Konstantin Preslavski", Universitetska str. 115, 9700 Shumen, Bulgaria}
\email{ivo\_michailov@yahoo.com, ziapkov2000@yahoo.co.uk}
\date{\today}
\thanks{This work is partially supported by a project of Shumen University}
\begin{document}
\baselineskip 20pt
\begin{abstract}
In this article we survey and examine the realizability of
$p$-groups as Galois groups over arbitrary fields. In particular we
consider various cohomological criteria that lead to necessary and
sufficient conditions for the realizability of such a group as a
Galois group, the embedding problem (i.e., realizability over a
given subextension), descriptions of such extensions, automatic
realizations among $p$-groups, and related topics.
\end{abstract}

\maketitle
\newcommand{\Gal}{{\rm Gal}}
\newcommand{\hw}{{\rm hw}}
\newcommand{\Ker}{{\rm Ker}}
\newcommand{\GL}{{\rm GL}}
\newcommand{\Br}{{\rm Br}}
\newcommand{\lcm}{{\rm lcm}}
\newcommand{\ord}{{\rm ord}}
\newcommand{\res}{{\rm res}}
\newcommand{\cor}{{\rm cor}}
\newcommand{\ch}{{\rm char}}
\newcommand{\tr}{{\rm tr}}
\newcommand{\Hom}{{\rm Hom}}
\newcommand{\ind}{{\rm ind}}
\newcommand{\Pin}{{\rm Pin}}
\newcommand{\Spin}{{\rm Spin}}
\newcommand{\attop}[2]{\genfrac{}{}{0pt}{1}{#1}{#2}}
\renewcommand{\thefootnote}{\fnsymbol{footnote}}
\numberwithin{equation}{section}

\section{Introduction: The inverse problem in Galois theory}
\label{1}

The purpose of this article is to survey and examine the
realizability of $p$-groups as Galois groups over arbitrary fields
for any prime $p$. We discuss the inverse problem of Galois theory,
its close relative -- the embedding problem, and related topics.

Let $G$ be a finite group, and let $K$ be a field. The \emph{inverse
problem of Galois theory} consists of two parts:
\begin{enumerate}
    \item {\bf Existence.} Determine whether there
    exists a Galois extension $M/K$ such that the Galois group $\Gal(M/K)$ is isomorphic to $G$.
    \item {\bf Actual construction.} If $G$ is realizable as a Galois
    group over $K$, construct explicitly either Galois extensions or polynomials over
    $K$ having $G$ as a Galois group.
\end{enumerate}

The classical inverse problem of Galois theory is the existence
problem for the field $K=\mathds Q$ of rational numbers. The
question whether all finite groups can be realized over $\mathds Q$
is one of the most challenging problems in mathematics, and it is
still unsolved.

In the early nineteenth century, the following result was
established:

\newtheorem{t1.1}{Theorem}[section]
\begin{t1.1}
{\rm (Kronecker-Weber)} {\it Any finite abelian group $G$ occurs as
a Galois group over $\mathds Q$. Furthermore, $G$ can be realized as
the Galois group of a subfield of the cyclotomic field $\mathds
Q(\zeta)$, where $\zeta$ is an $n$th root of unity for some natural
number $n$.}
\end{t1.1}

The proof can be found in most books on class field theory.

The first systematic study of the Inverse Galois Problem started
with Hilbert in 1892. Hilbert used his Irreducibility Theorem to
establish the following result:

\newtheorem{t1.2}[t1.1]{Theorem}
\begin{t1.2}
For any $n\geq 1$, the symmetric group $S_n$ and the alternating
group $A_n$ occur as Galois groups over $\mathds Q$.
\end{t1.2}

The first explicit examples  of polynomials with the alternating
group $A_n$ as a Galois group were given by Schur \cite{Schur} in
1930.

The next important step was taken in 1937 by A. Scholz and H.
Reichard \cite{Scho,Re} who proved the following existence result:

\newtheorem{t1.3}[t1.1]{Theorem}
\begin{t1.3}
For an odd prime $p$, every finite $p$-group occurs as a Galois
group over $\mathds Q$.
\end{t1.3}

The final step concerning solvable groups was taken by Shafarevich
\cite{Sha}, although with a mistake relative to the prime $2$. In
the notes appended to his Collected papers, p. 752, Shafarevich
sketches a method to correct this. For a full correct proof, the
reader is referred to the book by Neukirch, Schmidt and Wingberg
\cite[Chapter IX]{NSW}.

\newtheorem{t1.4}[t1.1]{Theorem}
\begin{t1.4}
{\rm (Shafarevich)} Every solvable group occurs as a Galois group
over $\mathds Q$.
\end{t1.4}

Extensive surveys of recent developments regarding the classical
inverse problem can be found in monographs such as
\cite{JLY,MM,Se-92,Vo}. We will, however, concentrate on the inverse
problem for $p$-groups over arbitrary fields.

Our paper is organized as follows. In Section \ref{2} we discuss the
cohomological approach to the embedding problem developed in such
works as \cite{ILF,MZ,Mi-cohII}. In Section \ref{3} we present some
more specific criteria concerning central embedding problems with
cyclic $p$-kernel. In Sections \ref{5} and \ref{6} we discuss the
quadratic corestriction homomorphism and orthogonal representations
of Galois groups. There we also give proofs of some unpublished
results of Michailov. In Section \ref{7} we present some of the most
significant results concerning the realizability of $p$-groups as
Galois groups over arbitrary fields. Finally, in Section \ref{8} we
investigate automatic realizations among $p$-groups.

\section{The embedding problem}
\label{2}

Let $k$ be arbitrary field and let $G$ is a non simple group. Assume
that $A$ is a normal subgroup of $G$. Then the realizability of the
quotient group $F=G/A$ as a Galois group over $k$ is a necessary
condition for the realizability of $G$ over $k$. In this way arises
the next generalization of the inverse problem in Galois theory --
the embedding problem of fields.

Let $K/k$ be a Galois extension with Galois group $F$, and let
\begin{equation}\label{e2.1}
1 \longrightarrow A \longrightarrow
G\overset{\alpha}{\longrightarrow} F \longrightarrow 1,
\end{equation}
be a group extension, i.e., a short exact sequence. Solving {\it the
embedding problem} related to $K/k$ and \eqref{e2.1} consists of
determining whether or not there exists a Galois algebra (called
also a {\it weak} solution) or a Galois extension (called a {\it
proper} solution) $L$, such that $K$ is contained in $L$, $G$ is
isomorphic to $\Gal(L/k)$, and the homomorphism of restriction to
$K$ of the automorphisms from $G$ coincides with $\alpha$. We denote
the so formulated embedding problem by $(K/k, G, A)$. We call the
group $A$ the {\it kernel} of the embedding problem.

A well known criterion for solvability is obtained by using the
Galois group $\Omega_k$ of the algebraic separable closure $\bar k$
over $k$.

\newtheorem{t2.1}{Theorem}[section]
\begin{t2.1}\label{t2.1}
\cite[Theorem 1.15.1]{ILF} The embedding problem $(K/k, G, A)$ is
weakly solvable if and only if there exists a homomorphism $\delta:
\Omega_k\to G$, such that $\alpha\cdot\delta=\varphi$, where
$\varphi: \Omega_k\to F$ is the natural epimorphism. The embedding
problem is properly solvable if and only if among the homomorphisms
$\delta$, there exists an epimorphism.
\end{t2.1}

Given that the kernel $A$ of the embedding problem is abelian,
another well known criterion holds. We can define an $F$-module
structure on $A$ by $a^\rho=\bar{\rho}^{-1}a\bar{\rho}$
($\bar{\rho}$ is a pre-image of $\rho\in F$ in $G$).

Let us recall the definition of the inflation map $\inf_F^{\Omega_k}
: H^2(F,A)\longrightarrow H^2(\Omega_k,A)$. Denote by
$G\times_F\Omega_k$ the direct product with amalgamated quotient
group $F$, i.e., the subgroup of the direct product $G\times
\Omega_k$, containing only the elements $(y,\omega)$, such that
$\alpha(y)=\varphi(\omega)$ for $y\in G$ and $\omega\in \Omega_k$.
We have then the commutative diagram with exact rows:

\begin{equation*}
\begin{CD}
1 @>>> A @>>> G\times_F\Omega_k @>\beta>> \Omega_k @>>> 1\\
  @.   @|   @V\psi VV @V\varphi VV  @.\\
1 @>>> A @>>> G @>\alpha>> F @>>> 1
\end{CD}
\end{equation*}
\\
where $\beta(y,\omega)=\omega$ and $\psi(y,\omega)=y$. The first row
gives us the inflation of the second row.

\newtheorem{c2.2}[t2.1]{Corollary}
\begin{c2.2}
\cite[Theorem 13.3.2]{ILF} Let $A$ be an abelian group and let $c$
be the $2$-coclass of the group extension \eqref{e2.1} in
$H^2(F,A)$. Then the embedding problem $(K/k, G,$ $A)$ is weakly
solvable if and only if $\inf_F^{\Omega_k}(c)=0$
\end{c2.2}

Next, let $K$ contain a primitive root of unity of order equal to
the order of the kernel $A$. Then we can define the character group
$\hat{A}=\text{Hom}\,(A,K^\ast)$ and make it an $F$-module by $
^\rho\chi(a)=\chi(a^\rho)^{\rho^{-1}}$, for $\chi\in\hat{A}$, $a\in
A$, $\rho\in F$.

Let $\mathds Z[\hat A]$ be the free abelian group with generators
$e_\chi$ (for $\chi\in \hat A$). We make it an $F$-module by $^\rho
e_\chi=e_{^\rho\chi}$. Then there exists an exact sequence of
$F$-modules
\begin{equation}\label{e2.2}
0 \longrightarrow V \longrightarrow \mathds Z[\hat
A]\overset{\pi}{\longrightarrow} \hat A \longrightarrow 0,
\end{equation}
where $\pi$ is defined by $\pi(\sum_ik_ie_{\chi_i^{\
}})=\prod_i\chi_i^{k_i}$ where $k_i\in \mathds Z$.

We can clearly consider all $F$-modules as $\Omega_k$-modules. The
exact sequence \eqref{e2.2} then implies the exact sequence
\begin{equation*}
0 \longrightarrow A\cong\Hom(\hat A,\bar k^\times) \longrightarrow
\Hom(\mathds Z[\hat A],\bar k^\times)\longrightarrow\Hom(V,\bar
k^\times) \longrightarrow 0.
\end{equation*}
Since $H^1(\Omega_k,\Hom(\mathds Z[\hat A],\bar k^\times))=0$ (see
\cite[\S 3.13.3]{ILF}), we obtain the following exact sequence
\begin{equation}\label{e2.3}
0 \longrightarrow H^1(\Omega_k,\Hom(V,\bar
k^\times))\overset{\beta}{\longrightarrow}
H^2(\Omega_k,A)\overset{\gamma}{\longrightarrow}H^2(\Omega_k,\Hom(\mathds
Z[\hat A],\bar k^\times)).
\end{equation}

We call the element $\eta=\gamma \bar c$ {\it the (first)
obstruction}. The condition $\eta=0$ clearly is necessary for the
solvability of the embedding problem $(K/k,G,A)$. This is the
well-known {\it compatibility} condition found by Faddeev and Hasse.
In general it is not a sufficient condition for solvability. Indeed
if we assume that $\eta=0$, then there appears a second obstruction,
namely $\xi\in H^1(\Omega_k,\Hom(V,\bar k^\times))$ such that
$\beta(\xi)=\bar c$. Thus, in order to obtain a necessary and
sufficient condition we must have both $\eta=0$ and $\xi=0$. The
second obstruction is very hard to calculate explicitly, though.
That is why embedding problems for which $H^1(\Omega_k,\Hom(V,\bar
k^\times))=0$ are of special interest. This condition turns out to
be fulfilled in a number of cases.

Let us begin with the so called \emph{Brauer problem}. The embedding
problem $(K/k, G, A)$ is called \emph{Brauer} if $\hat A$ is a
trivial $F$-module. Then we have the well known.

\newtheorem{t2.3}[t2.1]{Theorem}
\begin{t2.3}\label{t2.3}
{\rm (\cite[Theorem 3.1]{MZ}),\cite[Theorem 3.1]{ILF}} The
compatibility condition for the Brauer problem $(K/k, G, A)$ is
necessary and sufficient for its weak solvability.
\end{t2.3}

Recently, Michailov generalized this result with the following.

\newtheorem{t2.4}[t2.1]{Theorem}
\begin{t2.4}\label{t2.4}
{\rm (\cite[Theorem 3.2]{Mi-cohII})} Let $A$ be an abelian group of
order $n$, let the field $K$ contain a primitive $n$th root of
unity, and let $m$ be an integer such that $m^2\equiv1\pmod{n}$.
Assume that $(K/k,G,A)$ is an embedding problem such that the action
of $F=G/A$ on $\hat A$ satisfies the following: for any $\rho\in F$
we have either $^\rho\chi=\chi^m$ for all $\chi\in\hat A$, or
$^\rho\chi=\chi$ for all $\chi\in\hat A$. Then the compatibility
condition is necessary and sufficient for the weak solvability of
the embedding problem $(K/k,G,A)$.
\end{t2.4}

We can easily obtain now some previously known results as
corollaries from the latter theorem.

\newtheorem{c2.5}[t2.1]{Corollary}
\begin{c2.5}\label{c2.5}
{\rm (\cite[Theorem 3.2]{MZ})} Let the kernel $A$ be abelian and let
$ ^\rho\chi=\chi^{\pm1}$ for all $\chi\in \hat{A},\rho\in F$. Then
the compatibility condition is necessary and sufficient for the weak
solvability of the embedding problem $(K/k,G,A)$.
\end{c2.5}

\newtheorem{c2.6}[t2.1]{Corollary}
\begin{c2.6}\label{c2.6}
{\rm (\cite[\S 3.4.1]{ILF},\cite[Corollary 3.3]{MZ})} The embedding
problem $(K/k,G,A)$ with a kernel $A$ isomorphic to the cyclic group
of order $4$ is weakly solvable if and only if the compatibility
condition is satisfied.
\end{c2.6}

Since one of the forms of the compatibility condition is that all
associated Brauer problems are solvable, the above results reduce
the considerations of the original embedding problem to certain
associated Brauer problems. (For the definition of associated
problems see \cite{MZ}, and for the reduction see \cite{Mi-cohII}).

\section{Cohomological criteria for solvability of embedding problems with cyclic kernel of order $p$}
\label{3}

Let $k$ be arbitrary field of characteristic not $p$, containing a
primitive $p$th root of unity $\zeta$, and put
$\mu_p=\langle\zeta\rangle$. Let $K$ be a Galois extension of $k$
with Galois group $F$. Consider the group extension
\begin{equation}\label{e3.1}
1\longrightarrow \langle\varepsilon\rangle\longrightarrow
G\longrightarrow F\longrightarrow 1,
\end{equation}
where $\varepsilon$ is a central element of order $p$ in $G$. We are
going to identify the groups $\langle\varepsilon\rangle$ and
$\mu_p$, since they are isomorphic as $F$-modules.

Assume that $c\in H^2(F,\mu_p)$ is the $2$-coclass corresponding to
the group extension \eqref{e3.1}. \emph{The obstruction} to the
embedding problem $(K/k,G,\mu_p)$ we call the image of $c$ under the
inflation map $\inf_F^{\Omega_k}:H^2(F,\mu_p)\to
H^2(\Omega_k,\mu_p)$.

Note that we have the standard isomorphism of $H^2(\Omega_k,\mu_p)$
with the $p$-torsion in the Brauer group of $k$ induced by applying
$H^\ast(\Omega_k,\cdot)$ to the $p$-th power exact sequence of
$\Omega_k$-modules $1\longrightarrow\mu_p\longrightarrow\bar
k^\times\longrightarrow\bar k^\times\longrightarrow 1$. In this way,
the obstruction equals the equivalence class of the crossed product
algebra $(F,K/k,\bar c)$ for any $\bar c\in c$. Hence we may
identify the obstruction with a Brauer class in $\Br_p(k)$.

Note that we have an injection $\mu_p\hookrightarrow K^\times$,
which induces a homomorphism $\nu:H^2(F,\mu_p)\to H^2(F,K^\times)$.
Then the obstruction is equal to $\nu(c)$, since there is an
isomorphism between the relative Brauer group ${\rm Br}(K/k)$ and
the group $H^2(F,K^\times)$.

Clearly, the problem $(K/k,G,\mu_p)$ is Brauer, so from the proof of
Theorem \ref{t2.3} given in the paper \cite{MZ} it follows that
$H^1(\Omega_k,\Hom(V,\bar k^\times))=0$. Hence the homomorphism
$\gamma:H^2(\Omega_k,A)\to H^2(\Omega_k,\Hom(\mathds Z[\hat A],\bar
k^\times))$ is an injection. Therefore, the problem is solvable if
and only if the (first) obstruction is split.

More generally, the following result holds.

\newtheorem{t3.1}{Theorem}[section]
\begin{t3.1}\label{t3.1}
\cite{Ki} Let $c$ be the $2$-coclass in $H^2(F,\mu_p)$,
corresponding to the group extension \eqref{e3.1}. Then the
embedding problem $(K/k,G,\mu_p)$ is properly solvable if and only
if $\nu(c)=1$. If $K(\root p\of\beta)/k$ is a solution to the
embedding problem for some $\beta\in K^\times$, then all solutions
are of the kind $K(\root p\of{f\beta})/k$, for $f\in k^\times$.
\end{t3.1}

Henceforth, embedding problems of the kind $(K/k,G,\mu_p)$ we will
call for short $\mu_p$-embedding problems. An abstract description
of the solutions to the $\mu_p$-embedding problems is given by
Swallow in \cite{Sw1} and in a more concise form in \cite[Theorem
4]{ST}.

From the well-known Merkurjev-Suslin Theorem \cite{MeS} it follows
that the obstruction to any  $\mu_p$-embedding problem is equal to a
product of classes of $p$-cyclic algebras. The explicit computation
of these $p$-cyclic algebras, however, is not a trivial task. We are
going to discuss the methods for achieving this goal.

We denote by $(a, b;\zeta)$ the equivalence class of the $p$-cyclic
algebra which is generated by $i_1$ and $i_2$, such that $i_1^p=b,
i_2^p=a$ and $i_1i_2=\zeta i_2i_1$. For $p=2$ we have the quaternion
class $(a,b;-1)$, commonly denoted by $(a,b)$.

In 1987 Massy \cite{Ma-87} obtained a formula for the decomposition
of the obstruction in the case when $F=\Gal(K/k)$ is isomorphic to
$(C_p)^n$, the elementary abelian $p$-group.

\newtheorem{t3.2}[t3.1]{Theorem}
\begin{t3.2}\label{t3.2}
{\rm (\cite[Th\'eor\`eme 2]{Ma-87})} Let $K/k=k(\root p\of
{a_1},\root p\of {a_2},\dots,\root p\of {a_n})/k$ be a $(C_p)^n$
extension, and let $\sigma_1,\sigma_2,\dots,\sigma_n\in {\rm
Gal}(K/k)$ be given by $\sigma_i(\root p\of {a_j})/\root p\of
{a_j}=\zeta^{\delta_{ij}}$ ($\delta_{ij}$ is the Kronecker delta).
Let
\begin{equation*}
1\longrightarrow \mu_p\longrightarrow G\longrightarrow {\rm
Gal}(K/k)\longrightarrow 1
\end{equation*}
be a non split central extension, and choose pre-images
$s_1,s_2,\dots,s_n\in G$ of $\sigma_1,\sigma_2,\dots,\sigma_n$.
Define $d_{ij} (i\leq j)$ by $s_i^p=\zeta^{d_{ii}}$ and
$s_is_j=\zeta^{d_{ij}}s_js_i (i<j)$. Then the obstruction to the
embedding problem $(K/k,G,\mu_p)$ is
\begin{equation*}
\prod_{i=1}^{n}(a_i,\zeta;\zeta)^{d_{ii}}\prod_{i<k}(a_i,a_k;\zeta)^{d_{ik}}.
\end{equation*}
\end{t3.2}

In 1996 Swallow \cite{Sw2} described explicitly the solutions of
$\mu_p$-embedding problems with $F\cong (C_p)^n$.

In 2001 Quer \cite{Qu} found necessary and sufficient conditions for
solvability of the $\mu_p$-embedding problem with an abelian
quotient group $F$. The conditions are in terms of the existence of
elements with certain norm properties. They appear in the theory of
elliptic $\mathds Q$-curves discussed in the same paper.

In 2007 Michailov \cite{Mi-p4} obtained a formula for the
decomposition of the obstruction in the case when $F$ has a direct
factor $C_p$ for an odd $p$. The same is done in \cite{Mi-mod} for
$p=2$. We are now going to formulate these results in a unified way.

Let $H$ be a $p$-group and let
\begin{equation}\label{e3.2}
1\longrightarrow \mu_p\longrightarrow
G\overset{\pi}{\longrightarrow} F\cong H\times C_p\longrightarrow 1
\end{equation}
be a non split central group extension with characteristic
$2$-coclass $\gamma\in H^2(H\times C_p,C_p)$. By $\res_H\gamma$ we
denote the $2$-coclass of the group extension
\begin{equation*}
1\longrightarrow \mu_p\longrightarrow
\pi^{-1}(H)\overset{\pi}{\longrightarrow} H\longrightarrow 1.
\end{equation*}
Let $\sigma_1,\sigma_2,\dots,\sigma_m$ be a minimal generating set
for the maximal elementary abelian quotient group of $H$; and let
$\tau$ be the generator of the direct factor $C_p$. Finally, let
$s_1,s_2,\dots,s_m,t\in G$ be the pre-images of
$\sigma_1,\sigma_2,\dots,\sigma_m,\tau$, such that $t^p=\zeta^j$ and
$ts_i=\zeta^{d_i}s_it$, where $i\in\{1,2,\dots,m\};
j,d_i\in\{0,1,\dots,p-1\}$.

\newtheorem{t3.3}[t3.1]{Theorem}
\begin{t3.3}\label{t3.3}
{\rm(\cite[Theorem 4.1]{Mi-mod},\cite[Theorem 2.1]{Mi-p4})} Let
$K/k$ be a Galois extension with Galois group $H$ and let
$L/k=K(\root p\of b)/k$ be a Galois extension with Galois group
$H\times C_p$ ($b\in k^\times\setminus k^{\times p}$). Choose
$a_1,a_2,\dots,a_m\in k^\times$, such that $\sigma_k\root p\of
{a_i}=\zeta^{\delta_{ik}}\root p\of {a_i}$ ($\delta_{ik}$ is the
Kronecker delta). Then the obstruction to the embedding problem
$(L/k,G,\mu_p)$ is
\begin{equation*}
[K,H,\res_H\gamma]\left(b,\zeta^j\prod_{i=1}^{m}a_i^{d_i};\zeta\right).
\end{equation*}
\end{t3.3}

Ledet describes in his book \cite{Le-05} a more general formula for
the decomposition of the obstruction of $\mu_p$-embedding problems
with finite group $F$ isomorphic to a direct product of two groups.

Let $G$ be arbitrary finite group, and let $p$ be a prime divisor of
$\ord(G)$. Define $\mathcal O^p(G)$ as the subgroup of $G$ generated
by all elements of order prime to $p$. It is clear that $\mathcal
O^p(G)$ is the intersection of all normal subgroups in $G$ of
$p$-power index.

\newtheorem{t3.4}[t3.1]{Theorem}
\begin{t3.4}\label{t3.4}
{\rm(\cite[Theorem 6.1.4]{Le-05}} Let $L/k$ be a $N\times H$
extension, where $N$ and $H$ are finite groups. Let
\begin{equation*}
1\longrightarrow \mu_p\longrightarrow G\longrightarrow N\times
H\longrightarrow 1
\end{equation*}
be a non split central group extension with cohomology class
$\gamma\in H^2(N\times H,\mu_p)$. Let $K'/k$ and $K/k$ be the
subextensions corresponding to the factors $N$ and $H$. (I.e.,
$K'=L^H,K=L^N$.) Let $\sigma_1,\sigma_2,\dots,\sigma_m$ and
$\tau_1,\tau_2,\dots,\tau_n$ represent minimal generating sets for
the groups $N/\mathcal O^p(N)$ and $H/\mathcal O^p(H)$, and choose
$a_1,a_2,\dots,a_m,b_1,\dots b_n\in k^\times$, such that $\root p\of
{a_i}\in K^\times,\root p\of {b_i}\in K'^\times, \sigma_\kappa\root
p\of {a_i}=\zeta^{\delta_{i\kappa}}\root p\of {a_i}$ and
$\tau_\ell\root p\of {b_i}=\zeta^{\delta_{i\ell}}\root p\of {b_i}$
($\delta$ is the Kronecker delta). Finally, let
$s_1,\dots,s_m;t_1,\dots,t_n\in G$ be the pre-images of
$\sigma_1,\dots,\sigma_m;\tau_1,\dots,\tau_n$, and let
$d_{ij}\in\{0,\dots,p-1\}$ be given by
$t_js_i=\zeta^{d_{ij}}s_it_j$.

Then the obstruction to the embedding problem $(L/k,G,\mu_p)$ given
by $\gamma$ is
\begin{equation*}
[K,N,\res_N\gamma]\cdot [K',H,\res_H\gamma]\cdot
\prod_{i,j}(b_j,a_i;\zeta)^{d_{ij}}.
\end{equation*}
\end{t3.4}

There are many groups, however, which do not fit in the conditions
of the results listed above, so more criteria are needed. The
corestriction homomorphism $\cor_{F/H}:H^\ast (H,\cdot)\to H^\ast
(F,\cdot)$ appears to be one of the strongest tools in the case when
$F$ is not a direct product of smaller groups ($H$ is a properly
chosen subgroup of $F$). In \cite{Ri} is given an analogue of the
corestriction homomorphism, acting on central simple algebras. Given
a Galois extension $K_1/k$ of degree $p$, we have the corestriction
homomorhpism between the Brauer groups
$\cor_{K_1/k}:\Br(K_1)\longrightarrow \Br(k)$. Further, we denote by
$\Omega_k$ the Galois group of the separable closure of $k$ over the
field $k$, and by $\Omega_{K_1}$ the subgroup of $\Omega_k$, leaving
$K_1$ fixed. Then we also have the corestriction
$\cor_{\Omega_k/\Omega_{K_1}}:H^2(\Omega_{K_1},\mu_p)\longrightarrow
H^2(\Omega_k,\mu_p)$. Riehm shows in \cite[Theorem 11]{Ri} that the
following commutative diagram holds:

\begin{equation*}
\begin{CD}
H^2(\Omega_{K_1},\mu_p) @= \Br_p(K_1)\\
  @VV\cor_{\Omega_k/\Omega_{K_1}}V   @VV\cor_{K_1/k}V   \\
H^2(\Omega_k,\mu_p) @= \Br_p(k)
\end{CD}
\end{equation*}
\medskip

In this way one can easily see that the obstruction to the
corestricted embedding problem is equal to the corestriction of the
obstruction to the original problem. In order to apply the
corestriction homomorphism in the calculations, we need a formula
for the corestriction of $p$-cyclic algebras. Tignol gave in
\cite{Ti} a detailed proof of the so-called {\it projection
formula}: for $b\in k^\times\setminus k^{\times p}$ and $\delta\in
K_1$ we have
$\cor_{K_1/k}(\delta,b;\zeta)_{K_1}=(N_{K_1/k}(\delta),b;\zeta)_k$.
There are some cases, however, when the projection formula is not
enough. Swallow and Thiem found in \cite{ST} the following general
formula for the quadratic corestriction homomorphism.

\newtheorem{p}[t3.1]{Proposition}
\begin{p}\label{p}
{\rm(\cite[Proposition 4]{ST}} Let $a\in k^\times, K=k(\sqrt
a),\alpha_0=a_0+b_0\sqrt a$ and $\alpha_1=a_1+b_1\sqrt a$
($a_i,b_i\in k$). Then
\begin{enumerate}
\item If $b_{1-i}=0$, then
$\cor_{\Omega_k/\Omega_K}(\alpha_0,\alpha_1)_K=(a_{1-i},a_i^2-ab_i^2)_k$;
\item If $a_{1-i}b_i-a_ib_{1-i}=0$, then
$\cor_{\Omega_k/\Omega_K}(\alpha_0,\alpha_1)_K=(-a_ia_{1-i},a_i^2-ab_i^2)_k$;
\item Otherwise,\\
$\cor_{\Omega_k/\Omega_K}(\alpha_0,\alpha_1)_K=(a_0^2-ab_0^2,b_0(a_1b_0-a_0b_1))_k(a_1^2-ab_1^2,b_1(a_0b_1-a_1b_0))_k$.
\end{enumerate}
\end{p}

Next, we are going to describe a situation which sometimes occurs,
when a certain group extension can be constructed from other group
extension by lowering or raising the order of one of its generators.

Now, let $G$ be a finite group, and let
$\{\sigma_1,\dots,\sigma_k\}$ be a fixed (not necessarily minimal)
generating set of $G$ with these properties:
$\vert\sigma_1\vert=p^{n-1}$ for $n>1$, the subgroup $H$ generated
by $\sigma_2,\dots,\sigma_k$ is normal in $G$, and the quotient
group $G/H$ is isomorphic to the cyclic group $C_{p^{n-1}}$, i.e.,
$\sigma_1^i\notin H, 1\leq i<p^{n-1}$. Take now two arbitrary group
extensions
\begin{equation}\label{e3.3}
1 \longrightarrow \mu_p \longrightarrow G_1
\overset{\varphi}{\longrightarrow} G \longrightarrow 1
\end{equation}
and
\begin{equation}\label{e3.4}
1 \longrightarrow \mu_p\longrightarrow G_2
\overset{\psi}{\longrightarrow} G \longrightarrow 1.
\end{equation}
Denote by $\widetilde \sigma_i=\varphi^{-1}(\sigma_i)$ any preimage
of $\sigma_i$ in $G_1$ and by $\bar \sigma_i=\psi^{-1}(\sigma_i)$
any preimage of $\sigma_i$ in $G_2$, $i=1,\dots k$.

We write $G_2=G_1^{(p^n,\sigma_1)}$, if
\begin{enumerate}
\item $\vert\widetilde\sigma_1\vert=p^{n-1}$; \item
$\bar\sigma_1^{p^{n-1}}\in\mu_p,\bar\sigma_1^{p^{n-1}}\ne 1$; and
\item all other relations between the generators of the groups
$G_1$ and $G_2$ are identical, i.e.,
$\widetilde\sigma_i^{\alpha_i}=\zeta^l\prod_{j\ne
1}\widetilde\sigma_j^{\beta_j}\iff
\bar\sigma_i^{\alpha_i}=\zeta^l\prod_{j\ne 1}\bar\sigma_j^{\beta_j}$
 for $i=2,3,\dots,k$; $l,\alpha_i,\beta_j\in\mathds Z$; and
$[\widetilde\sigma_i,\widetilde\sigma_j]=\zeta^l\prod_{s\ne
1}\widetilde\sigma_s^{\varepsilon_s}\iff
[\bar\sigma_i,\bar\sigma_j]=\zeta^l\prod_{s\ne
1}\bar\sigma_s^{\varepsilon_s}$ for
$i,j=1,2,\dots,k;l,\varepsilon_s\in\mathds Z$.
\end{enumerate}

The following theorem, proved by Michailov \cite{Mi-coh}, gives us
the connection between the obstructions of the two embedding
problems related to \eqref{e3.3} and \eqref{e3.4}.

\newtheorem{t3.5}[t3.1]{Theorem}
\begin{t3.5}\label{t3.5}
{\rm(\cite[Theorem 2.7]{Mi-coh}} Let $L/F$ be a finite Galois
extension with Galois group $G=\Gal(L/F)$ as described above, let
$K=L^H$ be the fixed subfield of $H$, and let the groups $G_1$ and
$G_2$ from \eqref{e3.3} and \eqref{e3.4} be such that
$G_2=G_1^{(p^n,\sigma_1)}$. Denote by $O_{G_1}\in\Br_p(F)$ -- the
obstruction of the embedding problem $(L/F,G_1,\mu_p)$, by
$O_{G_2}\in\Br_p(F)$ -- the obstruction of the embedding problem
$(L/F,G_2,\mu_p)$, and by $O_{C_{p^n}}\in\Br_p(F)$ -- the
obstruction of the embedding problem $(K/F,C_{p^n},\mu_p)$ given by
the group extension $1 \longrightarrow \mu_p \longrightarrow C_{p^n}
\longrightarrow G/H\cong C_{p^{n-1}} \longrightarrow 1$. Then the
relation between these obstructions is given by
\begin{equation*}
O_{G_2}=O_{G_1}O_{C_{p^n}}\in\Br_p(F).
\end{equation*}
\end{t3.5}

In the following two sections we will focus on some specific
criteria related to the applications of the quadratic corestriction
homomorphism for $\mu_2$-embedding problems.

\section{The quadratic corestriction homomorphism}
\label{5}

Let $G$ be a finite group and let $H$ be an index $2$ subgroup of
$G$. Then we can define the quadratic corestriction homomorphism
$\cor_{G/H}:H^\ast(H,\cdot)\to H^\ast(G,\cdot)$. One way of
computing the quadratic corestriction is by applying the formula of
Tate \cite{Ta}, which in our case takes the form displayed in the
following

\newtheorem{l5.1}{Lemma}[section]
\begin{l5.1}\label{l5.1}
{\rm(\cite{Ta})} Let $G$ be a finite group, let $H$ be a subgroup of
index $2$ in $G$, let $g\in G\setminus H$ and let $\bar f\in
Z^2(H,\mu_2)$ be arbitrary $2$-cocycle. Define a map $f:G\times
G\longrightarrow \mu_2$ by:
\begin{equation*}
f(s_1,s_2)=\left\{
\begin{array}{ll}
\bar f(s_1,s_2)\bar f(gs_1g^{-1},gs_2g^{-1}), & if\ (s_1,s_2)\in H\times H\\
\bar f(s_1g^{-1},gs_2g^{-1})\bar f(gs_1,s_2), & if\ (s_1,s_2)\in
Hg\times H\\
\bar f(s_1,s_2g^{-1})\bar f(gs_1g^{-1},gs_2), & if\ (s_1,s_2)\in
H\times Hg\\
\bar f(s_1g^{-1},gs_2)\bar f(gs_1,s_2g^{-1}), & if\ (s_1,s_2)\in
Hg\times Hg.\\
\end{array}\right.
\end{equation*}
Then $f\in Z^2(G,\mu_2)$ and $[f]=\cor_{G/H}([\bar f])$, where
$[f]\in H^2(G,\mu_2)$ and $[\bar f]\in H^2(H,\mu_2)$ are  the
$2$-coclasses of $f$ and $\bar f$, respectively.
\end{l5.1}

With the aid of Lemma \ref{l5.1} we can construct some corestricted
group extensions and thus solve the related $\mu_2$-embedding
problems which can not be treated with the cohomological criteria
given in Section \ref{3}. However, Lemma \ref{l5.1} does not provide
us with the answer of the important question whether a given group
extension is a corestricted group extension, i.e., lies in the image
of $\cor_{G/H}$ for some subgroup $H$ of $G$. In what follows, we
will describe a construction which will help us relatively easy to
recognize the corestricted group extensions.

Let $\mathcal G$ be a finite $2$-group and let $E_4$ be a normal
subgroup of $\mathcal G$, isomorphic to the elementary abelian group
of order $4$ with generators $\sigma$ and $\tau$. Assume,
furthermore, that there exists a subgroup $\mathcal H$ in $\mathcal
G$, such that $E_4$ is a normal subgroup in $\mathcal H$, $\mathcal
H$ is contained in the centralizer $C_{\mathcal G}(E_4)$ of $E_4$ in
$\mathcal G$, and the index of $\mathcal H$ in $\mathcal G$ is $2$.
Next, choose and fix $g_1\in \mathcal G\setminus \mathcal H$, and
assume that $g_1\sigma g_1^{-1}=\sigma$ and $g_1\tau
g_1^{-1}=\sigma\tau$. Then for $H=\mathcal H/E_4$ and $G=\mathcal
G/E_4$ we have the isomorphism $G/H\cong \mathcal G/\mathcal H$.
Finally, choose and fix $g\in G\setminus H$, so that we have a
$G$-action on $E_4$, given by $c^h=c$ for all $c\in E_4$ and $h\in
H$; $\sigma^g=\sigma$ and $\tau^g=\sigma\tau$. In these notations,
it holds

\newtheorem{t5.2}[l5.1]{Theorem}
\begin{t5.2}\label{t5.2}
{\rm(\cite[Theorem 3.8]{Mi-coh})} Let $c_1\in H^2(G,\mu_2)$ be the
$2$-coclass, represented by the group extension $1\longrightarrow
E_4/\langle\sigma\rangle\cong \mu_2\longrightarrow \mathcal
G/\langle\sigma\rangle\longrightarrow G\longrightarrow 1$, let
$c_2\in H^2(H,\mu_2)$ be the $2$-coclass, represented by the group
extension $1\longrightarrow E_4/\langle\tau\rangle\cong
\mu_2\longrightarrow \mathcal H/\langle\tau\rangle\longrightarrow
H\longrightarrow 1$, and let $c_3\in H^2(H,\mu_2)$ be the
$2$-coclass, represented by the group extension $1\longrightarrow
E_4/\langle\sigma\tau\rangle\cong \mu_2\longrightarrow \mathcal
H/\langle\sigma\tau\rangle\longrightarrow H\longrightarrow 1$. Then
$\cor_{G/H}(c_2)=\cor_{G/H}(c_3)=c_1$.
\end{t5.2}

In order to verify whether a given group extension is corestricted
from some other group extension, we can try to construct the group
$\mathcal G$, having the properties given above Theorem \ref{t5.2}.
If we obtain defining relations that are not contradictory, then we
have reached our goal. This is done for instance in \cite[Section
5]{Mi-coh}. On the other hand, if we are not able to do that, it is
highly probable that the group extension is not a corestriction. In
that case we can use Proposition \ref{p5.4}, given below, in order
to prove such a supposition.

In the statements of the following results we use the standard
notations: $\ord(g)$ is the order of $g\in G$; $\exp(G)$ is the
exponent of the group $G$, i.e., $\exp(G)=\lcm\{\ord(g):g\in G\}$.
Clearly, when $G$ is a $2$-group, the exponent is equal to the
highest order of an element from $G$.

\newtheorem{l5.3}[l5.1]{Lemma}
\begin{l5.3}\label{l5.3}
Let
\begin{equation*}
\begin{CD} 1\longrightarrow \mu_2=\{\pm1\}\longrightarrow
Y @>\alpha>> Z \longrightarrow 1,
\end{CD}
\end{equation*}
be a group extension of $2$-groups, represented by the $2$-cocycle
$f\in Z^2(Z,\mu_2)$. Let $z\in Z, \ord(z)=k>1$ and let $y\in Y$ be
arbitrary preimage of $z$ in $Y$. Then
\begin{equation*}
f(z^{k/2},z^{k/2})=\left\{
\begin{array}{ll}
1, & if\ \ord(y)=k\\
-1, & if\ \ord(y)=2k.\\
\end{array}\right.
\end{equation*}
\end{l5.3}
\begin{proof}
As we know, the $2$-cocycle $f\in Z^2(Z,\mu_2)$ is defined by a set
$\{u_x\}_{x\in Z}\subset Y$ in this way: $u_xu_y=u_{xy}f(x,y)$ with
the natural condition $u_1=1,f(1,*)=f(*,1)=1$. Let $y$ and $z$ be
such as in the statement. We may assume that $y=u_z$. Clearly
$y^{k/2}=\pm u_{z^{k/2}}$, so
$y^k=y^{k/2}y^{k/2}=u_{z^{k/2}}u_{z^{k/2}}=u_{z^k}f(z^{k/2},z^{k/2})=f(z^{k/2},z^{k/2})$.
Since $\alpha(y^k)=z^k=1$, $y^k=\pm 1$ and we are done.
\end{proof}

\newtheorem{p5.4}[l5.1]{Proposition}
\begin{p5.4}\label{p5.4}
Let $G$ be a $2$-group, let $H$ be a subgroup of index $2$ in $G$
and let $g\in G\setminus H$. Further, let $\bar f\in Z^2(H,\mu_2)$
correspond to the group extension
\begin{equation*}
\begin{CD} 1\longrightarrow \mu_2=\{\pm1\}\longrightarrow
H_2 @>\beta>> H \longrightarrow 1,
\end{CD}
\end{equation*}
and let $f\in Z^2(G,\mu_2)$ correspond to the group extension
\begin{equation*}
\begin{CD} 1\longrightarrow \mu_2=\{\pm1\}\longrightarrow
G_1 @>\alpha>> G \longrightarrow 1,
\end{CD}
\end{equation*}
such that $[f]=\cor_{G/H}([\bar f])$. Put $H_1=\alpha^{-1}(H)$. Then
the following conditions hold:
\begin{enumerate}
\item $\exp(H_1)\leq\exp(H_2)$;
\item If for all $h\in H$ from
$\ord(h)=\exp(H)$ follows that $ghg^{-1}\in\langle h\rangle$, then
$\exp(H_1)=\exp(H)$.
\end{enumerate}
\end{p5.4}
\begin{proof}
$(1)$. If $\exp(H_1)=\exp(H)$, the inequality clearly holds. Now,
let $\exp(H_1)>\exp(H)$. Denote $k_1=\exp(H_1)$ and $k=\exp(H)$.
Then $k_1=2k$ and there exists $h_1\in H_1$ such that
$\ord(h_1)=k_1$. Also, for $h=\alpha(h_1)$ we have $\ord(h)=k$. From
Lemma \ref{l5.3} follows that $f(h^{k/2},h^{k/2})=-1$.

If $\exp(H_2)>\exp(H)$, then $\exp(H_2)=\exp(H_1)=k_1$. Now, suppose
$\exp(H_2)=\exp(H)=k$. Then there exists $h_2\in H_2$ such that
$\beta(h_2)=h$ and $\ord(h_2)=k$. Therefore $\bar
f(h^{k/2},h^{k/2})=1$. Next, consider the element $ghg^{-1}$, which
obviously has the same order $k$, and arbitrary preimage
$h_2'=\beta^{-1}(ghg^{-1})$. Then $\ord(h_2')=k=\exp(H_2)$ and from
Lemma \ref{l5.3} follows that $\bar
f(gh^{k/2}g^{-1},gh^{k/2}g^{-1})=1$. Finally, from Lemma \ref{l5.1}
we obtain $f(h^{k/2},h^{k/2})=\bar f(h^{k/2},h^{k/2})\bar
f(gh^{k/2}g^{-1},gh^{k/2}g^{-1})=1$, which is a contradiction.

$(2)$. Let $h\in H$ be of order $k=\exp(H)$. Then $ghg^{-1}=h^l$ for
some odd $l$. From Lemma \ref{l5.3} follows that $\bar
f(h^{k/2},h^{k/2})=\bar f((h^l)^{k/2},(h^l)^{k/2})=\bar
f(gh^{k/2}g^{-1},gh^{k/2}g^{-1})$, so $f(h^{k/2},h^{k/2})=1$.
Therefore any preimage of $h$ in $H_1$ will have order $k$. We are
done.
\end{proof}

The non abelian $2$-groups having a cyclic subgroup of index $2$ are
frequently discussed in various researches devoted to Galois theory.
For $n>3$, there are four such groups of order $2^n$ with exactness
up to an isomorphism (see for example \cite{Ha}): $Q_{2^n}$ (the
quaternion group), $D_{2^n}$ (the dihedral group)\footnote{We prefer
to denote by $D_{2m}$ the dihedral group of order $2m$, rather than
order $m$, to avoid any confusion with the notations for the
remaining groups.}, $SD_{2^n}$ (the semidihedral or quasidihedral
group) and $M_{2^n}$ (the modular group). These groups are generated
by two elements $\sigma$ and $\tau$. We list their presentations:
{\allowdisplaybreaks\begin{eqnarray*}
D_{2^n}&\cong& \langle\sigma,\tau\vert\sigma^{2^{n-1}}=\tau^2=1,\tau\sigma=\sigma^{-1}\tau\rangle\\
SD_{2^n}&\cong& \langle\sigma,\tau\vert\sigma^{2^{n-1}}=\tau^2=1,\tau\sigma=\sigma^{2^{n-2}-1}\tau\rangle\\
Q_{2^n}&\cong&
\langle\sigma,\tau\vert\sigma^{2^{n-1}}=1,\tau^2=\sigma^{2^{n-2}},\tau\sigma=\sigma^{-1}\tau\rangle\\
M_{2^n}&\cong&\langle \sigma, \tau\vert
\sigma^{2^{n-1}}=\tau^2=1,\tau\sigma=\sigma^{2^{n-2}+1}\tau\rangle.
\end{eqnarray*}}
If we suppose that the ground field contains enough roots of unity
(e.g. $2^{n-2}$-th roots of unity), the obstructions to the
$\mu_2$-embedding problems related to these groups are precisely
calculated (see \cite{Fr} for the dihedral and quaternion groups,
and \cite{Mi-ort} for the modular group). On the other hand, if we
relax the condition on the roots of unity even a little bit, the
calculations become very difficult. Michailov calculated the
obstructions to certain embedding problems with cyclic kernel in
\cite{Mi-dih,Mi-mod} under some specific requirements for the roots
of unity. If we have no requirements for the roots of unity it is
not known how to decompose the obstructions as products of
quaternion algebras for $n\geq 6$. That is why it is important to
know whether these groups can be constructed via a corestriction.
Unfortunately, as we will see, the answer is negative.

Note that the centre of $Q_{2^n}, D_{2^n},SD_{2^n}$ is
$\langle\sigma^{2^{n-2}}\rangle$, and the centre of $M_{2^n}$ is
$\langle\sigma^2\rangle$. Observe that if $G_1\cong Q_{2^n},
D_{2^n}$ or $SD_{2^n}$ then $G\cong D_{2^{n-1}}$, and if $G_1\cong
M_{2^n}$, then $G\cong C_{2^{n-2}}\times C_2$. In particular we see
that $\exp(G_1)=2^{n-1}$ and $\exp(G)=2^{n-2}$.

\newtheorem{p5.5}[l5.1]{Proposition}
\begin{p5.5}\label{p5.5}
The non abelian groups of order $2^n$ $(n>3)$, having a cyclic
subgroup of order $2^{n-1}$, can not be constructed via a
corestriction homomorphism.
\end{p5.5}
\begin{proof}
Assume that $G_1$ is one of these groups of order $2^n$ that can be
construted via a corestriction. Namely, let
\begin{equation*}
\begin{CD} 1\longrightarrow \langle\sigma^{2^{n-2}}\rangle\cong\mu_2\longrightarrow
G_1 @>\alpha>> G \longrightarrow 1,
\end{CD}
\end{equation*}
be a group extension corresponding to a $2$-cocylce $f\in
H^2(G,\mu_2)$ such that $[f]=\cor_{G/H}([\bar f])$ for a subgroup
$H<G$ of index $2$. Write as usual $H_1=\alpha^{-1}(H)$. There are
several possibilities for $H_1$ and $H=\alpha(H_1)$:
$H_1=\langle\sigma\rangle\cong C_{2^{n-1}}$ and
$H=H_1/\langle\sigma^{2^{n-2}}\rangle \cong C_{2^{n-2}}$; or
$H_1=\langle\sigma^2,\tau\rangle\cong D_{2^{n-1}}, Q_{2^{n-1}},
C_{2^{n-2}}\times C_2$ and $H=H_1/\langle\sigma^{2^{n-2}}\rangle
\cong D_{2^{n-2}}, C_{2^{n-3}}\times C_2$. In all possibilities,
however, the inequality $\exp(H_1)>\exp(H)$ holds.

Now, if $H\cong C_{2^{n-2}}$ or $D_{2^{n-2}}$ and $h\in H$ is an
element of maximal order, then all elements of maximal orders are in
$\langle h\rangle$. Indeed, in the dihedral groups we have the
relations $(\sigma^s\tau)^2=1$ for all $s$, so the elements of
maximal orders are only of this kind: $\sigma^l$ for an odd $l$.
From Proposition \ref{p5.4} (2) follows then that
$\exp(H_1)=\exp(H)$, a contradiction.

Finally, if $H\cong C_{2^{n-3}}\times C_2$ then $G\cong
C_{2^{n-2}}\times C_2$ and all elements from $G$ commute with the
elements from $H$. Then form Proposition \ref{p5.4} (2) again
follows that  $\exp(H_1)=\exp(H)$, which again is a contradiction.
\end{proof}

\section{Orthogonal representations of Galois groups}
\label{6}

We begin with some preliminaries about orthogonal representations.
Let $k$ be a field of characteristic $\ne 2$, let $V$ be a
finite-dimensional $k$-vector space, and let $(V,q)$ be a quadratic
space, $q$ being a quadratic form. The isometries $(V,q)\mapsto
(V,q)$ constitute a subgroup $O(q)$ of $\GL_k(V)$, called {\it the
orthogonal group} of $q$. {\it An orthogonal representation} of a
finite group $G$ is then a homomorphism $\mu:G\longrightarrow O(q)$
of $G$ into the orthogonal group of some regular quadratic form $q$.
From now on, by an orthogonal representation we will mean a {\it
faithful} one, i.e., an embedding $\mu: G\hookrightarrow O(q)$.

We adopt the notations about Clifford algebras used in \cite[Ch. 5,
S. 2]{Le-05}: $C(q)$ is the Clifford algebra of $q$; $C_0(q)$ is the
even Clifford algebra; $C(q)=C_0(q)\oplus C_1(q)$; if $x\in C_i(q)$,
we write $\partial x=i$; $C^\times(q)$ is the Clifford group,
defined as the subgroup of $C(q)^\times$, consisting of those
invertible elements $x$, for which $xVx^{-1}=V$. The anisotropic
vectors of $V$ are in $C^\times(q)$ and $vuv^{-1}=-T_v(u)$ for $u,
v\in V$, where $v$ is anisotropic and $T_v$ is the reflection on the
hyperplane $v^\perp$. There is an exact sequence
\begin{equation*}
\begin{CD} 1 \longrightarrow k^\times \longrightarrow C^\times(q) @>r>> O(q) \longrightarrow 1,
\end{CD}
\end{equation*}
$r$ being a map defined by $r_x:u\mapsto (-1)^{\partial x}xux^{-1}$,
where $x\in C^\times(q)$ and $u\in V$. In particular, for
$C_0^\times(q)=C^\times(q)\cap C_0(q)$ we get another exact sequence
\begin{equation*}
\begin{CD} 1 \longrightarrow k^\times \longrightarrow C_0^\times(q) @>r>> SO(q) \longrightarrow
1.
\end{CD}
\end{equation*}
Denote by $\iota$ the principal involution on $C(q)$, which
preserves the scalars, sums and vectors, and reverses products.
Denote by $N:C^\times(q)\longrightarrow k^\times$ the norm given by
$N(x)=x\iota(x)$, and by $sp:O(q)\longrightarrow k^\times/2$ the
spinor norm given by $sp(T_v)=\overline{q(v)}$. Put
$\Pin(q)=\ker(N), \Spin(q)=\Pin(q)\cap C_0^\times(q)$.

Hence, we have the long exact sequences
\begin{equation*}
\begin{CD} 1 \longrightarrow \mu_2 \longrightarrow \Pin(q) @>r>> O(q) @>sp>>
k^\times/2
\end{CD}
\end{equation*}
and
\begin{equation*}
\begin{CD} 1 \longrightarrow \mu_2 \longrightarrow \Spin(q) @>r>> SO(q) @>sp>>
k^\times/2.
\end{CD}
\end{equation*}
If we take the separable closure $\bar k$ of $k$, we get the short
exact sequences
\begin{equation}\label{8.1}
\begin{CD} 1 \longrightarrow \mu_2 \longrightarrow \Pin(\bar q) @>r>> O(\bar q)
\longrightarrow 1
\end{CD}
\end{equation}
and
\begin{equation}\label{8.2}
\begin{CD} 1 \longrightarrow \mu_2 \longrightarrow \Spin(\bar q) @>r>> SO(\bar q)
\longrightarrow 1.
\end{CD}
\end{equation}

\newtheorem{p6.1}{Proposition}[section]
\begin{p6.1}\label{p6.1}
Let $G$ and $H\leq G$ be as in the definitions above Theorem
\ref{t3.5} for $p=2$. Let $\varphi: G/H\longrightarrow \bar
k^\times$ be the homomorphism induced by the isomorphism
$G/H\cong\langle\zeta_{2^{n-1}}\rangle$ and the inclusion
$\langle\zeta_{2^{n-1}}\rangle\hookrightarrow\bar k^\times$, where
$\zeta_{2^{n-1}}$ is a primitive $2^{n-1}$th root of unity. Assume
that is given an orthogonal representation $G\hookrightarrow O(q)$,
and take the restriction
$1\longrightarrow\mu_2\longrightarrow\widetilde G\longrightarrow
G\longrightarrow 1$ of
$1\longrightarrow\mu_2\longrightarrow\Pin(\bar q)\longrightarrow
O(\bar q)\longrightarrow 1$. Then there exists a subgroup $\bar G$
of $C^\times(\bar q)$, such that
\begin{enumerate}
    \item The diagram\\
\begin{equation*}
\begin{CD}
1 @>>> \mu_2 @>>> \bar G @>r'>> G @>>> 1\\
  @.   @VVNV   @VVNV @VV\varphi V  @.\\
1 @>>> 1 @>>> \bar k^\times @= \bar k^\times @>>> 1
\end{CD}
\end{equation*}
\\
is commutative with exact rows, where $N$ is the norm and $r'$ is
the restriction of $r:C^\times(\bar q)\longrightarrow O(\bar q)$ on
$\bar G$;
    \item Either $\bar G=\widetilde G^{(2^n,\sigma_1)}$, or $\widetilde G=\bar G^{(2^n,\sigma_1)}$.
\end{enumerate}
\end{p6.1}
\begin{proof}
Choose and fix preimages
$\tilde\sigma_1,\dots,\tilde\sigma_\kappa\in \widetilde G$ of the
generators $\sigma_1,\dots,\sigma_\kappa$ of $G$. Next, let $\bar G$
be the subgroup of $C^\times(\bar q)$, generated by the elements
$\bar\sigma_1=\tilde\sigma_1\zeta_{2^n},\bar\sigma_2=\tilde\sigma_2,\dots,\bar\sigma_\kappa=\tilde\sigma_\kappa$.
We will show that $\bar G$ satisfies the conditions (1) and (2).

First, we will see that $\ker(r')\cong\mu_2$. Choose arbitrary $x$
from $\ker(r')$, i.e., $x=\prod
\bar\sigma_1^{i_1}\bar\sigma_2^{i_2}\cdots\bar\sigma_\kappa^{i_\kappa}$
and $r'(x)=\prod
\sigma_1^{i_1}\sigma_2^{i_2}\cdots\sigma_\kappa^{i_\kappa}=1$. Since
$\sigma_1^{i_1}\notin H$ for $1\leq i_1<2^{n-1}$, we get $i_1\equiv
0\ (\text{mod}\ 2^{n-1})$. Whence
$x=\pm\bar\sigma_2^{i_2}\cdots\bar\sigma_\kappa^{i_\kappa}=\pm\tilde\sigma_2^{i_2}\cdots\tilde\sigma_\kappa^{i_\kappa}\in\ker(r)\cong\mu_2$,
where by $r$ we denote also the restriction of $r$ on $\Pin(\bar
q)$.

Next,
$N(\bar\sigma_1)=N(\tilde\sigma_1)\zeta_{2^n}^2=\zeta_{2^{n-1}}=\varphi
r'(\bar\sigma_1)$ and $N(\bar\sigma_i)=\varphi r'(\bar\sigma_i)=1$
for $i=2,\dots,k$. Therefore, the diagram in the statement indeed is
commutative and with exact rows.

Finally, we have either $\vert\tilde\sigma_1\vert=2^{n-1}$, or
$\tilde\sigma_1^{2^{n-1}}=-1$. If
$\vert\tilde\sigma_1\vert=2^{n-1}$, we have
$\bar\sigma_1^{2^{n-1}}=\zeta_{2^n}^{2^{n-1}}=-1$, and since the
remaining relations between the generators of $\widetilde G$,
respectively of $\bar G$, are identical, we see that $\bar
G=\widetilde G^{(2^n,\sigma_1)}$. If $\tilde\sigma_1^{2^{n-1}}=-1$,
we have  $\bar\sigma_1^{2^{n-1}}=1$, so $\widetilde G=\bar
G^{(2^n,\sigma_1)}$.
\end{proof}

We will look now at the double covers of the symmetric group $S_n$.
Let $L/k$ be a Galois extension with Galois group $G$, and assume
that $L$ is the splitting field over $k$ of an irreducible
polynomial $f(x)\in k[x]$ of degree $n$. We can then embed $G$
transitively into $S_n$ by considering the elements of $G$ as
permutations of the roots of $f(x)$. Consider the 'positive' double
cover
\begin{equation*}
1\longrightarrow \mu_2\longrightarrow \widetilde
S_n^+\longrightarrow S_n\longrightarrow 1,
\end{equation*}
where transpositions lift to elements of order $2$, and the
'negative' double cover
\begin{equation*}
1\longrightarrow \mu_2\longrightarrow \widetilde
S_n^-\longrightarrow S_n\longrightarrow 1,
\end{equation*}
where the transpositions lift to elements of order $4$ (in both
cases products of two disjoint transpositions lift to elements of
order $4$). Take now the restrictions
\begin{equation*}
1\longrightarrow \mu_2\longrightarrow \widetilde G^+\longrightarrow
G\longrightarrow 1,
\end{equation*}
and
\begin{equation*}
1\longrightarrow \mu_2\longrightarrow \widetilde G^-\longrightarrow
G\longrightarrow 1,
\end{equation*}
of the 'positive' and, respectively, the 'negative' double covers of
$S_n$.

\newtheorem{c6.2}[p6.1]{Corollary}
\begin{c6.2}
Under the above notations, assume that $G$ contains a transposition.
Then $\widetilde G^-=\widetilde G^{+(4,\sigma_1)}$ and the relation
between the obstructions of the related embedding problems is given
by $O_{\widetilde G^-}=(-1,d_f)O_{\widetilde G^+}$, where $d_f$ is
the discriminant of $f(x)$.
\end{c6.2}
\begin{proof}
Denote by $\sigma_1$ any transposition from $G$, and by $H$ the
subgroup $G\cap A_n$ of $G$. Choose a set of generators
$\sigma_2,\dots,\sigma_\kappa$ for $H$, so
$\sigma_1,\sigma_2,\dots,\sigma_\kappa$ are generators for $G$. Then
pick their preimages
$\tilde\sigma_1,\tilde\sigma_2,\dots,\tilde\sigma_k$ in $\widetilde
G^+$. Since $\widetilde G^+$ is in $\Pin_n(\bar k)$ (see
\cite{Le-05}) and $\sigma_1^2=\tilde\sigma_1^2=1$, we can make use
of Proposition \ref{p6.1}. Whence we can define a subgroup $\bar
G=\widetilde G^{+(4,\sigma_1)}$ of $C_n^\times(\bar k)$, generated
by $\bar\sigma_1=\tilde\sigma_1 i,
\bar\sigma_2=\tilde\sigma_2,\dots,\bar\sigma_\kappa=\tilde\sigma_\kappa$
($i$ is the imaginary unity). Now, from $\bar\sigma_1^2=-1$ follows
that each transposition from $G$ lifts to an element of order $4$ in
$\bar G$. Indeed, if $\sigma$ is another transposition in $G$, then
$\sigma=\sigma_1\tau$ for $\tau\in H$, so we can choose preimage
$\bar\sigma=\bar\sigma_1\bar\tau=i\tilde\sigma_1\bar\tau\in \bar G$,
where $\bar\tau$ is also in $\widetilde G^+$, and
$\ord(\tilde\sigma_1\bar\tau)=2$. Therefore, $\bar\sigma^2=-1$ and
$\bar G\cong\widetilde G^-$.

Finally, from $\sigma_1(\sqrt{d_f})=-\sqrt{d_f}$ we obtain the
relation between the obstructions, given in the statement.
\end{proof}

Now, let us recall the definition of Galois twist, which involves
the existence of the first cohomological group $H^1(G,\mathcal G)$,
where $\mathcal G$ is non abelian group with a $G$-action. Assume
again that $(V,q)$ is a quadratic space over $k$, and that $K/k$ is
a Galois extension with Galois group $G$. Then we can extend the
scalars to get a quadratic space $(V_K,q_K)$. The semi-linear action
of $G$ then gives us the equation $q_K(\sigma u)=\sigma q_K(u)$.
Conversely, if $(W,Q)$ is a quadratic space over $K$ endowed with a
semi-linear action such that $Q(\sigma u)=\sigma Q(u)$ is satisfied,
we obtain a quadratic space $(W^G,Q^G)$ over $k$ by taking fixed
points and restricting $Q$. These two operations (scalar extension
and fixed points) preserve regularity and are each others inverses.
Also, $O(Q)$ is a $G$-group by conjugation:
$(\sigma\varphi)(u)=\sigma\varphi(\sigma^{-1}u)$.

Next, let $f:G\to O(q_K)$ be a crossed homomorphism. Then we can
define a semi-linear action by $^\sigma u=f_\sigma(\sigma u)$ and
get an induced quadratic space $(V_f,q_f)=((V_K)^G,(q_K)^G)$ over
$k$. Furthermore, if $g$ is equivalent to $f$, i.e.,
$g_\sigma=\varphi f_\sigma\sigma\varphi^{-1}$ for some $\varphi\in
O(q_K)$, then $V_g=\varphi(V_f)$, and consequently $(V_f,q_f)$ and
$(V_g,q_g)$ are equivalent. Hence, to each element in $H^1(G,O(q))$
we can associate an equivalence class of quadratic spaces over $k$.

The quadratic space $(V_f,q_f)$ is said to arise from $(V,q)$ by
taking the {\it Galois twist} with respect to $f$.

Define the element
\begin{equation*}
\hw(q)=\prod_{i<j}(a_i,a_j)\in\Br(k),
\end{equation*}
where $a_i=q(u_i)$ for some canonical orthogonal basis
$u_1,\dots,u_n$ of $q$. Clearly, $\hw(q)$ depends only on the
equivalence class of $q$. It is called the {\it Hasse-Witt
invariant} or the {\it second Stiefel-Whitney class} of $q$.

The obstruction now can be calculated by the formula, displayed in
the following.

\newtheorem{t6.3}[p6.1]{Theorem}
\begin{t6.3}\label{t6.3}
{\rm (\cite{Fr,Le-05})} Let $L/k$ be a finite Galois extension with
Galois group $G=\Gal(L/k)$ and assume $G\hookrightarrow O(q)$ for
some regular quadratic form $q$ over $k$. Let $e:\Gal(\overline
k/k)\to O(q)$ be the induced crossed homomorphism, and let $q_e$ be
the Galois twist of $q$ by $e$. Also, let
\begin{equation*}
1\longrightarrow \mu_2\longrightarrow \widetilde G\longrightarrow
G\longrightarrow 1
\end{equation*}
be the group extension induced by $G\hookrightarrow O(q)$ and the
group extension
\begin{equation*}
\begin{CD} 1 \longrightarrow \mu_2 \longrightarrow \Pin(\overline q) @>r>> O(\overline q)
\longrightarrow 1.
\end{CD}
\end{equation*}
Let $K/k=k(\sqrt{a_1},\dots,\sqrt{a_r})/k$ be the maximal elementary
abelian $2$-subextension of $L/k$, and let $\rho_1,\dots,\rho_r\in
G$ be such that
$\rho_i(\sqrt{a_j})=(-1)^{\delta_{ij}}\cdot\sqrt{a_j}$. Then the
obstruction to the embedding problem $(L/k,\widetilde G,\mu_2)$ is
\begin{equation*}
\hw(q)\hw(q_e)(d,-d_e)\prod_{i=1}^r(a_i,sp(\rho_i))\in\Br(k),
\end{equation*}
where $d$ and $d_e$ are discriminants of $q$ and $q_e$, repectively.
\end{t6.3}

Now, let $L/k$ be a finite Galois extension with Galois group $G$,
let $H$ be a subgroup of $G$ with fixed field $K=L^H$, and let
$\mu:H\hookrightarrow O(q)$ be an orthogonal representation over
$k$. Then, according to \cite{Fr,FM}, we can construct an {\it
induced orthogonal representation} $\ind\mu: G\hookrightarrow
O(q_{\ind\mu})$, where $\ind\mu$ has as underlying module the
induced $G$-module of the $H$-module $V_q:V_{\ind\mu}=\oplus
(V_q\otimes\sigma)=V_q\otimes_{kH}kG$, $\sigma$ running over a given
right transversal $R$ of $H$ in $G$. Note that $V_q\subset
V_{\ind\mu}$ is a subspace which is $H$-invariant. It is not hard to
show that, given an orthogonal representation $\mu:H\hookrightarrow
O(q)$, such $V_{\ind\mu}$ exists and is unique up to an isomorphism
(see e.g. \cite[\S 3.3]{FH}). Moreover, the action of $G$ can be
explicitly determined: Each element $v\in V_{\ind\mu}$ has a unique
expression $v=\sum w_\sigma\otimes\sigma$ for elements $w_\sigma$ in
$V_q$. For a given $g\in G$, we must have
\begin{equation}\label{e1.1}
g\cdot(w_\sigma\otimes\sigma)=hw_\sigma\otimes\tau\ \ \text{if}\
g\sigma=\tau h\ (\tau\in R).
\end{equation}

Next, assume that we have a special orthogonal representation
$\mu:H\hookrightarrow SO(q)$ over $k$. Denote by $\bar k$ the
separable closure of $k$, and by $\bar q$ the extension of $q$ to
$\bar k$. Then we have a diagram
\begin{equation*}
\begin{CD}
1 @>>> \mu_2 @>>> \widetilde H @>>> H @>>> 1\\
  @.   @VVV   @VVV @VVV  @.\\
1 @>>> \mu_2 @>>> \Spin(\bar q) @>>> SO(\bar q) @>>> 1\\
  @.   @VVV   @VVV @VVV  @.\\
1 @>>> \bar k^\times @>>> C_0^\times(\bar q) @>>> SO(\bar q) @>>> 1,
\end{CD}
\end{equation*}
\\
where as usual $\mu_2=\{\pm 1\}$ and
$1\longrightarrow\mu_2\longrightarrow\widetilde H\longrightarrow
H\longrightarrow 1$ is the restriction of
$1\longrightarrow\mu_2\longrightarrow\Spin(\bar q)\longrightarrow
SO(\bar q)\longrightarrow 1$. The induced orthogonal representation
$\ind\mu: G\hookrightarrow O(q_{\ind\mu})$, in its turn, gives us
the diagram
\begin{equation*}
\begin{CD}
1 @>>> \mu_2 @>>> \widetilde G @>>> G @>>> 1\\
  @.   @VVV   @VVV @VVV  @.\\
1 @>>> \mu_2 @>>> \Pin(\bar q_{\ind\mu}) @>>> O(\bar q_{\ind\mu}) @>>> 1\\
  @.   @VVV   @VVV @VVV  @.\\
1 @>>> \bar k^\times @>>> C^\times(\bar q_{\ind\mu}) @>>> O(\bar
q_{\ind\mu}) @>>> 1.
\end{CD}
\end{equation*}
\\
Recently, Michailov \cite{Mi-ort} proved the following.

\newtheorem{t6.4}[p6.1]{Theorem}
\begin{t6.4}\label{t6.4}
{\rm (\cite[Theorem 2.2]{Mi-ort})} Let $G$ be a finite group, and
let $H$ be a subgroup of $G$, such that $\vert H\vert=2^tm, (t,m\geq
1)$. Let also $\mu:H\hookrightarrow SO(q)$ be an orthogonal
representation over $k$ with an underlying module $V_q$, such that
$n=\dim_k V_q\equiv0\ ({\rm mod}\ 4)$. Denote by $\bar f\in
Z^2(H,\mu_2)$ and by $f\in Z^2(G,\mu_2)$ the $2$-cocycles given by
the described above group extensions
$1\longrightarrow\mu_2\longrightarrow\widetilde H\longrightarrow
H\longrightarrow 1$ and
$1\longrightarrow\mu_2\longrightarrow\widetilde G\longrightarrow
G\longrightarrow 1$, respectively. Then $[f]=\cor_{G/H}([\bar f])$,
where $\cor_{G/H}:H^2(H,\mu_2)\longrightarrow H^2(G,\mu_2)$ is the
corestriction map.
\end{t6.4}

Now, assume again that $L/k$ is a normal and separable extension
with a finite Galois group $G$. We can always find a primitive
element $\theta$ such that $L=k(\theta)$. Let $f(x)\in k[x]$ be the
minimal polynomial of $\theta$ of degree $n=[L:k]$, and let
$\theta=\theta_1,\theta_2,\dots,\theta_n$ be the conjugates of
$\theta$. Then $G=G(f)$ embeds transitively into the symmetric group
$S_n$.

For a given proper subgroup $H$ of $G$, we set $m=\vert H\vert$ and
$\kappa=(G:H)=n/m$. Clearly, $\theta$ is a primitive element of the
extension $L/K$ as well, where $K=L^H$. Since the minimal polynomial
of $\theta$ over $K$ divides $f(x)$, we can assume that
$\theta=\theta_1,\theta_2,\dots,\theta_m$ for $1<m=[L:K]<n$ are the
conjugates of $\theta$ over $K$. $H$ embeds transitively in $S_m$,
so we can take the group extension
\begin{equation}\label{e1.1.1}
1\longrightarrow \mu_2\longrightarrow \widetilde H\longrightarrow
H\longrightarrow 1,
\end{equation}
which is the restriction of the group extension
\begin{equation*}
1\longrightarrow \mu_2\longrightarrow \widetilde S_m\longrightarrow
S_m\longrightarrow 1,
\end{equation*}
$\widetilde S_m$ being the positive double cover of $S_m$.

Next, recall that for the quadratic form $q_1=\langle
1,\dots,1\rangle$ on $V_1=k^m$ we have that $S_m$ embeds in
$O_m(k)=O(q_1)$, so we get an orthogonal representation
$H\hookrightarrow O_m(k)$. Set $q=q_1\perp q_2\perp\dots\perp
q_\kappa$ and $V=V_1\oplus V_2\oplus\dots\oplus V_\kappa$, where
$q_1=q_2=\dots=q_\kappa$ and $V_1=V_2=\dots V_\kappa$. In this way,
we get the induced orthogonal representation $G\hookrightarrow
O_n(k)$, which is identical to the transitive embedding of $G=G(f)$
in $S_n$. Now, take the group extension
\begin{equation}\label{e1.1.2}
1\longrightarrow \mu_2\longrightarrow \widetilde G\longrightarrow
G\longrightarrow 1,
\end{equation}
which is the restriction of the group extension
\begin{equation*}
1\longrightarrow \mu_2\longrightarrow \widetilde S_n\longrightarrow
S_n\longrightarrow 1.
\end{equation*}
Denote by $\overline f\in Z^2(H,\mu_2)$ the $2$-cocycle representing
\eqref{e1.1.1} and by $f\in Z^2(G,\mu_2)$ the $2$-cocycle
representing \eqref{e1.1.2}, i.e., $\overline f=\res(s_m)$ and
$f=\res(s_n)$. From Theorem \ref{t6.4} now follows that
$[f]=\cor_{G/H}([\overline f])$, under the extra assumptions
$H\hookrightarrow SO_m(k)$ and $m\equiv0\ ({\rm mod}\ 4)$.

\section{Realizability of $p$-groups as Galois groups}
\label{7}

The realizability of $2$-groups as Galois groups is discussed in a
great number of papers, but obstructions over arbitrary base field
$k$ are obtained in only a few of them. The precise calculation of
the obstructions in general is a difficult and delicate task, but
when done properly, it rewards the author with the opportunity to
describe the Galois extensions realizing that group, and to discover
new automatic realizations (see Section \ref{8}).

The condition for the realizability of the cyclic group $C_4$ over
arbitrary fields with characteristic different from $2$ has been
known for a long time. Namely, if $k(\sqrt a)/k$ is an arbitrary
quadratic extension, then the group $C_2=\Gal(k(\sqrt a)/k)$ can be
embedded into a $C_4$ extension if and only if $a$ is a sum of two
squares. Hence $C_4$ is realizable over $k$ if and only if there
exists $a\in k^\times\setminus k^{\times2}$ such that $a=x^2+y^2$
for some $x,y\in k$. It is easy to see that this condition descends
from the obstruction to the embedding problem $(k(\sqrt
a)/k,C_4,\mu_2)$ which is exactly the quaternion algebra
$(a,-1)\in\Br_2(k)$.

In 1886 Dedekind constructed for the first time Galois extensions
over $\mathds Q$ with Galois group $Q_8$, the quaternion group of
order $8$. This result was published in a posthumous article
\cite{De}. In 1936 Witt published a paper \cite{Wi} where he
obtained a general description of quaternion extensions over
arbitrary fields with characteristic equal to or different from $2$.
Moreover, Witt found a necessary and sufficient condition for the
realizability of $p$-groups over fields with characteristic $p$. In
1984  J.-P. Serre generalized Witt's construction of quaternion
extensions \cite[n. 3.2, Remarque]{Se-84}. A year later appeared a
very significant paper by Fr\"ohlich \cite{Fr} where he developed
the theory of orthogonal representations of Galois groups,
generalizing Serre's formula from \cite{Se-84}. In particular,
Fr\"ohlich calculated the obstructions to realizability of the
dihedral and quaternion groups of order $4n$ over fields containing
a primitive $n$-th root of unity for $n\geq 2$ (see
\cite[(7.10)]{Fr}. An explicit description of the solutions to
embedding problems associated to orthogonal Galois representations
was done by Crespo in \cite{Cr}.

The first attempt at finding necessary and sufficient conditions for
the realizability of a number of non abelian groups of order $16$
was made by Kiming \cite{Ki} in 1990. Five years later Ledet
\cite{Le-95} found the decomposition of the obstructions as
quaternion algebras of all groups of order $2^n$ for $n\leq 4$ with
the exception of $C_{16}$, the cyclic group of order $16$. For some
of the groups Ledet applied Theorem \ref{t3.4}, and for others (e.g.
the dihedral, semi-dihedral and the quaternion groups) he applied
'brute force'. The latter approach exploits the well-known theorem
for the centralizer of a c.s. subalgebra:

{\it If $A$ is a c.s. $k$-algebra and $B$ is a c.s. $k$-subalgebra
of $A$ then $A\cong B\otimes_k C_A(B)$, where $C_A(B)$ is the
centralizer of $B$ in $A$}.

Unfortunately, this method works exclusively for Brauer problems
such that the obstruction can be decomposed as a product of at most
three quaternion algebras (i.e., when the quotient group is of order
$\leq 8$, since the crossed product algebra then has dimension $\leq
64=4^3$ over the base field). Ledet made a parametrization of all
Galois extensions realizing groups of order $16$ in \cite{Le-01} in
the case when the obstructions are decomposed as a product of two
quaternion algebras. Such a description can be done for any
$2$-group if its obstruction is of this kind (see \cite{Mi-32}
concerning one group of order $32$ and \cite{Mi-Q16} for some
particular cases concerning the quaternion group of order $16$). In
\cite{JLY} the reader can find a good survey of explicit
descriptions of generic polynomials and extensions realizing
$2$-groups as Galois groups.

Michailov \cite{Mi-dih,Mi-mod} applied brute force for specific
Brauer problems with cyclic $2$-kernel involving the four non
abelian $2$-groups having a cyclic subgroup of index $2$.

The obstruction to the embedding of a $C_8$ extension into a
$C_{16}$ extension was found by Swallow \cite{Sw-C16} by applying
brute force again. There are many more results concerning the
realizability of cyclic and abelian $2$-groups in articles such as
\cite{Schn,MSc,AFSS,JP}.

Grundman, Smith and Swallow wrote an extensive survey \cite{GSS} of
the results known until 1995 concerning the groups of order $2^n$
for $n\leq 4$.

Naturally, the next goal was the investigation of groups of order
$32$. There are $51$ groups of order $32$. With the aid of the
corestriction homomorphism Swallow and Thiem calculated in \cite{ST}
the obstructions to the realizability of several non abelian groups
of orders $32$ and $64$.

By applying a variation of Theorem \ref{t3.3}, Grundman and Stewart
\cite{GSt} found the obstructions to the realizability of $13$
groups of order $32$ that have a quotient group isomorphic to a
direct product of cyclic groups of order $\leq 4$ and/or the
dihedral group of order $8$. T. Smith \cite{Sm} investigated some
realizability properties (e.g. obstructions, solutions and automatic
realizations) of the two extra-special groups of order $32$.

A systematic study of all non abelian groups of order $32$ was done
by Michailov \cite{Mi-32} in 2007. Grundman and Smith further
elaborated on the obstructions of some of these groups in
\cite{GS-D4,GS-32}. The latter two authors also determined in
\cite{GS-64} the obstructions to the realizability of $134$ non
abelian groups of order $64$. They used a variation of Theorem
\ref{t3.4} to calculate the obstructions for $34$ $\mu_2$-embedding
problems with quotient $(C_2)^r\times (C_4)^s\times (D_8)^t$. The
remaining $100$ obstructions were obtained for groups known as
\emph{pullbacks} by a method applied by Michailov in \cite{Mi-32}
for $18$ groups of order $32$.\footnote{In fact, the first examples
of the obstructions to the realizability of pullbacks were given by
Ledet in \cite{Le-95}.}

Namely, let $\varphi' : G'\rightarrow F$ and $\varphi'' :
G''\rightarrow F$ be homomorphisms with kernels $N'$ and,
respectively, $N''$. {\it The pullback} of the pair of homomorphisms
$\varphi'$ and $\varphi''$ is called the subgroup in $G'\times G''$
of all pairs $(\sigma',\sigma'')$, such that
$\varphi'(\sigma')=\varphi''(\sigma'')$. The pullback is denoted by
$G'\curlywedge G''$. It is also called the direct product of the
groups $G'$ and $G''$ with amalgamated quotient group $F$ and
denoted by $G'*_F G''$.

Now, let $N_1=N'\times \{1\}$ and $N_2=\{1\}\times N''$. Then $N_1$
and $N_2$ are normal subgroups of $G'\curlywedge G''$, such that
$N_1\cap N_2=\{1\}$. The converse is also true (see \cite{ILF}, I,
\S 12):

\newtheorem{l4.1}{Lemma}[section]
\begin{l4.1}\label{l4.1}
Let $N_1$ and $N_2$ be two normal subgroups of the group $G$, such
that $N_1\cap N_2=\{1\}$. Then $G$ is isomorphic to the pullback
$(G/N_1)\curlywedge (G/N_2)$.
\end{l4.1}

The application to embedding problems is given by:

\newtheorem{t4.0}[l4.1]{Theorem}
\begin{t4.0}\label{t4.0}
{\rm (\cite[Theorem 1.12]{ILF})} Let $K/k$ be a Galois extension
with Galois group $F$. In the notations of the lemma, let $F\cong
G/N_1N_2$ and $G\cong (G/N_1)\curlywedge (G/N_2)$. Then the
embedding problem $(K/k, G, N_1\times N_2)$ is solvable if and only
if the embedding problems $(K/k, G/N_1, N_2)$ and $(K/k, G/N_2,
N_1)$ are solvable.
\end{t4.0}

Next, according to a paper by Ninomia \cite{Ni}, there are $26$ non
isomorphic non abelian groups of order $2^n$ for $n\geq 4$ that have
a cyclic subgroup of index $4$. Michailov determined in
\cite{Mi-coh} the obstructions to the realizability of $14$ groups
by applying Theorems \ref{t3.3}, \ref{t3.5} and \ref{t5.2}. The
obstructions to the remaining groups were obtained in
\cite{Mi-cohII} by taking embedding problems with cyclic kernel of
order $2^{n-3}$ and applying Theorem \ref{t2.4}. This is done with
the assumption that the base field $k$ contains a primitive
$2^{n-3}$th root of unity.
\medskip

In what follows we will concentrate on results about $p$-groups for
an odd prime $p$, although most of the groups have their twins for
$p=2$.

We begin with a $\mu_p$-embedding problem, involving the cyclic
group $C_{p^2}$. Let $K=k(\root p\of a)$, where $a\in
k^\times\setminus k^{\times p}$. Now, consider the embedding problem
given by $K/k$ and the group extension:
\begin{equation}\label{e4.1}
1\longrightarrow \mu_p\longrightarrow C_{p^2}\longrightarrow
C_p\longrightarrow 1.
\end{equation}
Denote by $\Gamma$ the crossed product algebra, corresponding to
\eqref{e4.1}. As we know, $\Gamma$ is generated by elements $\root
p\of a$ and $u$, such that $(\root p\of a)^p=a, u\root p\of
a=\zeta\root p\of au$ and $u^p=\zeta$. Therefore
$[\Gamma]=(a,\zeta;\zeta)$, so $[\Gamma]=1$ if and only if $\zeta\in
N_{K/k}(K^\times)$.

Naturally, the two non abelian groups of order $p^3$ were the first
non abelian $p$-groups investigated for realizability as Galois
groups. The first one is the Heisenberg group of exponent $p$. We
denote it by $G_1$ and its generators by $g_1,g_2$ and $g_3$, such
that $g_1^p=g_2^p=g_3^p=1,g_1g_2=g_2g_1g_3$ and $g_3$ is central.
The second group we denote it by $G_2$. It is generated by $g_1$ and
$g_2$, such that $g_1^{p^2}=g_2^p=1$ and $g_1g_2=g_2g_1^{p+1}$.

Massy \cite{Ma-87} investigated the realizability of $G_1$ and $G_2$
over arbitrary fields containing a primitive $p$-th root of unity.

\newtheorem{t4.1}[l4.1]{Theorem}
\begin{t4.1}\label{t4.1}
{\rm(\cite{Ma-87},\cite[Theorem 3.1]{Mi-p4})} Let $K=k(\root p\of
a_1,\root p\of a_2)$ be a $C_p\times C_p$ extension of $k$ and let
$K_i=k(\root p\of a_i), i=1,2$. Denote the generators of $C_p\times
C_p$ by $\sigma_1$ and $\sigma_2$, such that $\sigma_i(\root p\of
{a_j})/\root p\of {a_j}=\zeta^{\delta_{ij}}, i=1,2$. The embedding
problem given by $K/k$ and the group extension
\begin{equation*}
1\longrightarrow \langle g_3\rangle\cong \mu_p\longrightarrow
G_1\underset{\attop{g_1\mapsto\sigma_1}{
g_2\mapsto\sigma_2}}{\longrightarrow} C_p\times C_p\longrightarrow 1
\end{equation*}
is solvable if and only if $a_2\in N_{K_1/k}(K_1^\times)$. In that
case for $\omega\in K_1^\times$, such that $N(\omega)=a_2$ we put
$\alpha=\omega^{p-1}\sigma_1(\omega)^{p-2}\cdots\sigma_1^{p-2}(\omega)$.
Then the set $\{K(\root p\of {f\alpha})\mid f\in k^\times\}$ gives
all solutions.
\end{t4.1}

\newtheorem{t4.2}[l4.1]{Theorem}
\begin{t4.2}\label{t4.2}
{\rm(\cite{Ma-87},\cite[Theorem 3.2]{Mi-p4})} Let $K=k(\root p\of
a_1,\root p\of a_2)$ be a $C_p\times C_p$ extension of $k$ and let
$K_i=k(\root p\of a_i), i=1,2$. Denote the generators of $C_p\times
C_p$ by $\sigma_1$ and $\sigma_2$, such that $\sigma_i(\root p\of
{a_j})/\root p\of {a_j}=\zeta^{\delta_{ij}}, i=1,2$. The embedding
problem given by $K/k$ and the group extension
\begin{equation*}
1\longrightarrow \langle g_1^p\rangle\cong \mu_p\longrightarrow
G_2\underset{\attop{g_1\mapsto\sigma_1}{
g_2\mapsto\sigma_2}}{\longrightarrow} C_p\times C_p\longrightarrow 1
\end{equation*}
is solvable if and only if $a_2\zeta\in N_{K_1/k}(K_1^\times)$. In
that case for $\omega\in K_1^\times$, such that $N(\omega)=a_2\zeta$
we put
$\alpha=\omega^{p-1}\sigma_1(\omega)^{p-2}\cdots\sigma_1^{p-2}(\omega)$.
Then the set $\{K(\root p\of {f\root p\of {a_1}^{-1}\alpha})\mid
f\in k^\times\}$ gives all solutions.
\end{t4.2}

More results about the realizability of $G_1$ and $G_2$ over global
fields were obtained in \cite{MN}.

Michailov \cite{Mi-p4} investigated four non abelian groups of order
$p^4$, that have a quotient group of the kind $H\times C_p$.
Generally, there are $15$ groups of order $p^4$ for any odd prime
$p$. Of course, we do not need to bother about the groups of the
type $G\times H$, since their realizability depends only on the
realizability of $G$ and $H$. From these four groups three have as a
quotient group the group $C_{p^2}\times C_p$ and one has $(C_p)^3$
as a quotient group. The remaining non abelian groups of order $p^4$
do not have a quotient group that is a direct product of smaller
groups.

Firstly, let us describe the non abelian groups, having
$C_{p^2}\times C_p$ as a quotient group, which is generated by
$\sigma_1$ and $\sigma_2$ such that $\sigma_1^{p^2}=\sigma_2^p=1$
and $\sigma_1\sigma_2=\sigma_2\sigma_1$. The presentations of these
groups can be given by the relations between their generators
$g_1,g_2,g_3$ and $g_4$. The symbol $[a,b]$ below stands for the
commutator $a^{-1}b^{-1}ab$.
\begin{eqnarray*}
&&G_3 : g_1^p=g_4,g_2^p=g_3^p=g_4^p=1, [g_2,g_1]=g_3,\ g_3\ {\rm and}\ g_4\ {\rm are\ central},\\
&&G_4 : g_1^p=g_4, g_2^p=g_3, g_3^p=g_4^p=1, [g_2,g_1]=g_3,\ g_3\ {\rm and}\ g_4\ {\rm are\ central},\\
&&G_5 : g_1^p=g_3, g_3^p=g_4, g_2^p=g_4^p=1, [g_2,g_1]=g_4,\ g_3\
{\rm and}\ g_4\ {\rm are\ central}.
\end{eqnarray*}
The corresponding central group extensions are as follows:
\begin{eqnarray}
\label{4.1} &&1\longrightarrow \langle g_3\rangle\cong
\mu_p\longrightarrow G_3\underset{\attop{g_1\mapsto\sigma_1}{
g_2\mapsto\sigma_2}}{\longrightarrow} C_{p^2}\times
C_p\longrightarrow
1,\\
\label{4.2} &&1\longrightarrow \langle g_3\rangle\cong
\mu_p\longrightarrow G_4\underset{\attop{g_1\mapsto\sigma_1}{
g_2\mapsto\sigma_2}}{\longrightarrow} C_{p^2}\times
C_p\longrightarrow
1,\\
\label{4.3} &&1\longrightarrow \langle g_4\rangle\cong
\mu_p\longrightarrow G_5\underset{\attop{g_1\mapsto\sigma_1}{
g_2\mapsto\sigma_2}}{\longrightarrow} C_{p^2}\times
C_p\longrightarrow 1.
\end{eqnarray}

The last group is given by the presentation
\begin{equation*}
G_6 : g_1^p=g_2^p=1, g_3^p=g_4, g_4^p=1, [g_2,g_1]=g_4,\ g_3\ {\rm
and}\ g_4\ {\rm are\ central}.
\end{equation*}
Given that the group $(C_p)^3$ is generated by elements $\rho_1,
\rho_2$ and $\rho_3$, we obtain the central group extension:
\begin{equation}\label{4.4}
1\longrightarrow \langle g_4\rangle\cong \mu_p\longrightarrow
G_6\underset{g_i\mapsto\rho_i}{\longrightarrow}
(C_p)^3\longrightarrow 1.
\end{equation}

Let $k$ be a field with characteristic $\ne p$, let $\zeta$ be a
primitive $p$th root of unity in $k$, and let $a_1,a_2\in k^\times$
be linearly independent mod $k^{\times p}$. Denote $K=k(\root p\of
{a_1},\root p\of {a_2})$ and $K_i=k(\root p\of {a_i}), i=1,2$. Now,
assume that the embedding problem $(k(\root p\of {a_1})/k,
C_{p^2},\mu_p)$ is solvable. Then $(a_1,\zeta;\zeta)=1$, so there
exists $\alpha\in K_1$, such that $\zeta=N_{K_1/k}(\alpha)$. Let
$L_1/k=K_1(\root p\of {f_1\beta})/k$ be arbitrary $C_{p^2}$
extension, where $f_1\in k^\times$ and $\beta=\root p\of
{a_1}(\alpha^{p-1}\sigma_1(\alpha)^{p-2}\cdots$
$\sigma_1^{p-2}(\alpha))^{-1}$. Then we have a $C_{p^2}\times C_p$
extension $L=L_1(\root p\of {a_2})$. We display the obstructions and
the solutions to the $\mu_p$-embedding problems in the following
three theorems.

\newtheorem{t4.3}[l4.1]{Theorem}
\begin{t4.3}\label{t4.3}
{\rm(\cite[Theorem 4.1]{Mi-p4})} The obstruction to solvability of
the embedding problem given by the Galois extension $L/k$ and the
group extension \eqref{4.1} is $(a_2,a_1;\zeta)$. If the embedding
problem is solvable, i.e., $a_2=N_{K_1/k}(\omega)$ for $\omega\in
K_1^\times$, we may put
$\gamma=\omega^{p-1}\sigma_1(\omega)^{p-2}\cdots
\sigma_1^{p-2}(\omega)$. Then all Galois extensions realizing $G_3$
are $\{L(\root p\of {f_2\gamma})/k : f_2\in k^\times\}$.
\end{t4.3}

\newtheorem{t4.4}[l4.1]{Theorem}
\begin{t4.4}\label{t4.4}
{\rm(\cite[Theorem 4.2]{Mi-p4})} The obstruction to solvability of
the embedding problem given by the Galois extension $L/k$ and the
group extension \eqref{4.2} is $(a_2,a_1\zeta;\zeta)$. If the
embedding problem is solvable, i.e., $a_1\zeta=N_{K_2/k}(x)$ for
$x\in K_2^\times$, we may put $\omega=\root p\of
{a_2}(x^{p-1}\sigma_2(x^{p-2})\sigma_2^2(x^{p-3})\cdots
\sigma_2^{p-2}(x))^{-1}$. Then all Galois extensions realizing $G_4$
are $\{L(\root p\of {f_2\omega})/k : f_2\in k^\times\}$.
\end{t4.4}

\newtheorem{t4.5}[l4.1]{Theorem}
\begin{t4.5}\label{t4.5}
{\rm(\cite[Theorem 4.3]{Mi-p4})} The obstruction to solvability of
the embedding problem given by the Galois extension $L/k$ and the
group extension \eqref{4.3} is $[L_1, C_{p^2},
\zeta](a_2,a_1;\zeta)$. If a primitive $p^2$th root of unity
$\zeta_{p^2}=\root p\of \zeta$ is contained in $k$, then the
obstruction is $(\zeta_{p^2}^{-1}a_2, a_1;\zeta)$. Given that the
embedding problem is solvable, i.e.,
$\zeta_{p^2}^{-1}a_2=N_{K_1/k}(y)$, for some $y\in K_1$, we may put
$\omega=\root {p^2}\of {a_1}y^{p-1}\sigma_1(y)^{p-2}\cdots
\sigma_1^{p-2}(y)$. Then all Galois extensions realizing $G_5$ are
$\{L(\root p\of {f\omega})/k : f\in k^\times\}$.
\end{t4.5}

The groups of orders $p^5$ and $p^6$ that have abelian quotients
obtained by factoring out $\mu_p$ or $(\mu_p)^2$ will be
investigated in a future paper by Michailov \cite{Mi-p5}.

In 2009 Michailov investigated in \cite{Mi-ort} the modular
$p$-group $M(p^n)$ and the group extensions from
$H^2(M(p^n)),\mu_2)$, for $n\geq 3$. The modular $p$-group $M(p^n)$
is generated by two elements $\alpha$ and $\beta$, such that
$\alpha^{p^{n-1}}=\beta^p=1$ and
$\beta\alpha=\alpha^{1+p^{n-2}}\beta$, for $n\geq 3$. Put
$q=p^{n-2}$.

\newtheorem{t4.6}[l4.1]{Theorem}
\begin{t4.6}\label{t4.6}
{\rm(\cite[p. 3723]{Mi-ort})} Let $a_1,a_2\in k^\times$ be
independent mod $k^{\times p}$ and denote $K_i=k(\root p\of {a_i}),
i=1,2$. Let $K/k$ be a $C_q=\langle\sigma\rangle$ extension, such
that $K_1\subset K$. Then $L/k=K(\root p\of {a_2})/k$ is a
$C_q\times C_p$ extension, generated by elements $\sigma$ and
$\tau$, such that $\sigma^q=\tau^p=1$. Consider the group extension:
\begin{equation*}
1\longrightarrow \mu_p\cong\langle \alpha^q\rangle\longrightarrow
M(p^n)\underset{\attop{\alpha\mapsto\sigma}{\beta\mapsto\tau}}{\longrightarrow}
C_q\times C_p\longrightarrow 1.
\end{equation*}
The obstruction to solvability of the embedding problem
$(L/k,M(p^n),\mu_p)$ is
\begin{equation}\label{3.3}
[K,C_q,\zeta](a_2,a_1;\zeta)\in \Br(k),
\end{equation}
where $[K,C_q,\zeta]$ is the equivalence class of the crossed
product cyclic algebra $(K,\sigma,\zeta)$, given by the restricted
group extension
\begin{equation*}
1\longrightarrow \mu_p\cong\langle \alpha^q\rangle\longrightarrow
C_{pq}\underset{\attop{\alpha\mapsto\sigma}{}}{\longrightarrow}
C_q\longrightarrow 1.
\end{equation*}
\end{t4.6}

According to \cite[p. 3278]{Mi-ort}, the cohomological group
$H^2(M(p^n),\mu_p)$ is isomorphic to $\mu_p^2$, and consists of the
$2$-coclasses related to the group extensions
\begin{equation*}
1\longrightarrow \mu_p\longrightarrow
G_{\varepsilon_1,\varepsilon_2}\underset{\attop{\tilde\alpha\mapsto\alpha}{\tilde\beta\mapsto\beta}}{\longrightarrow}
M(p^n)\longrightarrow 1,
\end{equation*}
where
$\varepsilon_i\in\mu_p,\tilde\beta^p=\varepsilon_1,\tilde\alpha^q[\tilde\beta,\tilde\alpha]=\varepsilon_2$
and $\tilde\alpha^{pq}=1$. With the aid of the corestriction
homomorphism, Michailov determined in \cite{Mi-ort} the obstructions
to the $\mu_p$-embedding problems, related to these group
extensions.

\newtheorem{t4.7}[l4.1]{Theorem}
\begin{t4.7}\label{t4.7}
{\rm(\cite[Proposition 4.3]{Mi-ort})} Let $L/k$ be a $M(p^n)$
extension containing the biquadratic extension $k(\root p\of
{a_1},\root p\of {a_2})\subset L$, where $K_1=k(\root p\of
{a_1})=L^{C_q\times C_p}$. The obstruction to the embedding problem
$(L/k,G_{1,\zeta},\mu_p)$,  where $G_{1,\zeta}\cong \langle x,y,z :
x^{p^{n-1}}=y^p=z^p=1, z -\ \text{central}, yx=x^{q+1}yz\rangle$ is
$(a_2,a_1;\zeta)\in \Br_p(k)$.
\end{t4.7}

\newtheorem{t4.8}[l4.1]{Theorem}
\begin{t4.8}\label{t4.8}
{\rm(\cite[Proposition 4.4]{Mi-ort})} Let $L/k$ be a $M(p^n)$
extension containing the biquadratic extension $k(\root p\of
{a_1},\root p\of {a_2})\subset L$, where $K_1=k(\root p\of
{a_1})=L^{C_q\times C_p}$. The obstruction to the embedding problem
$(L/k,G_{\zeta,1},\mu_p)$, where $G_{\zeta,1}\cong \langle x,y :
x^{p^{n-1}}=y^{p^2}=1, y^p -\ \text{central}, yx=x^{q+1}y\rangle$ is
$(a_2,\zeta;\zeta)\in \Br_p(k)$.
\end{t4.8}

\newtheorem{t4.9}[l4.1]{Theorem}
\begin{t4.9}\label{t4.9}
{\rm(\cite[Proposition 4.5]{Mi-ort})} Let $L/k$ be a $M(p^n)$
extension containing the biquadratic extension $k(\root p\of
{a_1},\root p\of {a_2})\subset L$, where $K_1=k(\root p\of
{a_1})=L^{C_q\times C_p}$. The obstruction to the embedding problem
$(L/k,G_{\zeta,\zeta},\mu_p)$, where $G_{\zeta,\zeta}\cong \langle
x,y : x^{p^{n-1}}=y^{p^2}=1, y^p -\ \text{central},
yx=x^{q+1}y^{p+1}\rangle$ is $(\zeta a_1,a_2;\zeta)\in \Br_p(k)$.
\end{t4.9}

Michailov also made in \cite{Mi-ort} an explicit (to some extent)
description of the modular $p$-extensions.

\newtheorem{t4.10}[l4.1]{Theorem}
\begin{t4.10}\label{t4.10}
{\rm(\cite[Theorem 3.2]{Mi-ort})} Let $K_1=k(\root p\of {a_1})$ for
$a_1\in k^\times\setminus k^{\times p}$, and let $K/k$ be a cyclic
$C_q=\langle\sigma\rangle$ extension, such that $K_1\subset K$.
$L/k$ is an $M(p^n)$ Galois extension, solving the embedding problem
$(K/k,M(p^n),C_p\times C_p)$, if and only if, there exist $b_0\in
K^\times\setminus K^{\times p}, f\in k^\times\setminus k^\times\cap
K^{\times p}$ and $x\in K^\times$ such that $\sigma(b_0)/b_0=fx^p,
L/k=K(\root p\of {b_0},\root p\of f)/k$ and $c=f^{q/p}N_{K/k}(x)$ is
a $p$th root of unity, but $c\ne 1$
\end{t4.10}

Waterhause's ideas from \cite{Wa} were further developed in
\cite{Mi-ort} with a number of propositions, e.g. the following.

\newtheorem{p4.10}[l4.1]{Proposition}
\begin{p4.10}\label{p4.10}
{\rm(\cite[Proposition 3.4]{Mi-ort})} Let
$L=K(b_0^{1/p},b_1^{1/p},b_2^{1/p})$, where $b_0\in
K^\times\setminus K^{\times p}, b_1=\sigma(b_0)/b_0\in
K^\times\setminus K^{\times p}, b_2=\sigma(b_1)/b_1\in
K^\times\setminus K^{\times p}$ and $\sigma(b_2)/b_2\in K^{\times
p}$. Then $L/k$ is Galois and the Galois group of $L/k$ is
isomorphic either to the semidirect product  $(C_p)^3\rtimes C_q$,
or to the group, generated by $\sigma_1,\tau_0,\tau_1$ and $\tau_2$,
such that
$\sigma_1^q=\tau_0,\sigma_1\tau_1\sigma_1^{-1}=\tau_1\tau_0^{-1},\sigma_1\tau_2\sigma_1^{-1}=\tau_2\tau_0\tau_1^{-1}$,
where $\tau_0,\tau_1$ and $\tau_2$ are the generators of
$\Gal(L/K)$, given by $\tau_i(\root p\of {b_j})=\root p\of
{b_j}\zeta^{\delta_{ij}}$.
\end{p4.10}

Note that for $q=p$, the group $(C_p)^3\rtimes C_p$ from the latter
proposition is a non abelian group of order $p^4$ that is not
isomorphic to any of the groups $G_3,\dots,G_6$ considered so far.
Denote it by $G_7$. It is generated by elements
$\sigma,\tau,\lambda$ and $\mu$ with the following relations:
\begin{equation*}
G_7\cong (C_p)^3\rtimes C_p :
\sigma^p=\tau^p=\lambda^p=\mu^p=1,[\lambda,\tau]=1,[\mu,\tau]=\sigma,[\mu,\lambda]=\tau,\
\sigma\ {\rm is\ central}.
\end{equation*}

A necessary and sufficient condition for the embedding of a
$C_p$-extension into a $G_7$-extension will be obtained in the
example given at the end of this Section.

It turns out that the description of $p$-extensions given so far has
an impact on certain realizability issues over the field of rational
numbers $\mathds Q$. By applying Michailov's methods from
\cite{Mi-p4,Mi-ort}, S. Checcoli \cite{Ch} provided recently a
characterization of infinite algebraic Galois extensions of the
rationals with uniformly bounded local degrees.

\newtheorem{t4.11}[l4.1]{Theorem}
\begin{t4.11}\label{t4.11}
{\rm(\cite[Theorem 1]{Ch})} Let $K/\mathds Q$ be an infinite Galois
extension. Then the following conditions are equivalent:
\begin{enumerate}
    \item $K$ has uniformly bounded local degrees at every prime;
    \item $K$ has uniformly bounded local degrees at almost every prime;
    \item $\Gal(K/\mathds Q)$ has finite exponent.
\end{enumerate}
Moreover, if $K/\mathds Q$ is abelian, then the three properties:

$\ \ (a)\ K$ has uniformly bounded local degrees;

$\ \ (b)\ K$ is contained in $Q^{(d)}$ for some positive integer d;

$\ \ (c)\ $ every finite subextension of $K$ can be generated by
elements of bounded degree;\\ are equivalent. However, in general,
we have that (c) implies (b) which implies (a) and none of the
inverse implications holds.
\end{t4.11}

Waterhause's ideas \cite{Wa} also inspired J. Min\'a\v c and J.
Swallow, who considered in \cite{MS-05} certain embedding problems
with a cyclic $p$-kernel. We are now going to state their main
result.

Let $k$ be an arbitrary field, and suppose that $K/k$ is a cyclic
extension with Galois group $G=\Gal(K/k)\cong \mathds Z/p\mathds
Z$\footnote{Henceforth we hold on to the additive notation $\mathds
Z/n\mathds Z$ for the cyclic group of order $n$, which is used by
Min\'a\v c, Schultz and other authors.}, with generator $\sigma$.
Let $A=\oplus_{j=0}^{p-1}\mathds F_p\tau^j$ be a free $\mathds
F_p[G]$-module on the generator $\tau$, where $\sigma$ acts by
multiplication by $\tau$. Let $A_i$ be the $\mathds
F_p[G]$-submodule generated by $(\tau-1)^i$. Finally, let $\mathcal
E_i,1<i\leq p$, denote the Galois embedding problem related to the
group extension:
\begin{equation*}
\mathcal E_i:\quad 1\to A_1/A_i\to (A/A_i)\rtimes G\to
(A/A_1)\rtimes G=\Gal(L/k)\to 1.
\end{equation*}
Observe that $A/A_1\cong\mathds F_p$, a trivial $\mathds
F_p[G]$-module; hence $(A/A_1)\rtimes G\cong \mathds Z/p\mathds
Z\times \mathds Z/p\mathds Z$. We also assume that the projection of
$(A/A_1)\rtimes G=\Gal(L/k)$ onto $G$ coincides with the restriction
map $\Gal(L/k)\to G=\Gal(K/k)$. Assume that $\ch\ k\ne p$ and denote
by $\zeta$ any primitive $p$-th root of unity. Then $\Gal(L/k)$ is
naturally isomorphic to $\Gal(L(\zeta)/k(\zeta))$. After identifying
these two Galois groups we set $k(\zeta,\root p\of b)$ to be the
fixed field of $1\rtimes G$ in $L(\zeta)$.

\newtheorem{t4.12}[l4.1]{Theorem}
\begin{t4.12}\label{t4.12}
{\rm(\cite[Theorem 1]{MS-05})}

{\bf A)} Let $k$ be an arbitrary field. Then the following are
equivalent:
\begin{enumerate}
    \item Some $\mathcal E_i$ is solvable.
    \item Each $\mathcal E_i$ is solvable.
\end{enumerate}
Consequently, if $(A/A_2)\rtimes G$ occurs as a Galois group over
$k$, then $(A/A_i)\rtimes G$ occurs as well, for all $2\leq i\leq
p$.

{\bf B)} Now assume that $\ch\ k\ne p$. Then $(1)$ and $(2)$ are
also equivalent to

$\ \ (3)\ b\in N_{K(\zeta)/k(\zeta)}(K(\zeta)^\times)$.

{\bf C)} Now assume further that $\zeta\in k$. Suppose $(1)-(3)$
hold, and let $\omega\in K^\times$ satisfy $N(\omega)=b$. Suppose
$i>2$. Then a solution to $\mathcal E_i$ is given by
\begin{equation*}
\tilde L=K(\root p\of{f\omega^{(\sigma-1)^{p-i}}},\root
p\of{\omega^{(\sigma-1)^{p-i+1}}},\dots,\root
p\of{\omega^{(\sigma-1)^{p-2}}}),
\end{equation*}
$f\in k^\times$. If $i=2$ then a solution to $\mathcal E_2$ is given
by $\tilde L=K(\root p\of{\omega^{(\sigma-1)^{p-2}}})$. Moreover,
all solutions of $\mathcal E_i$ arise in this way.
\end{t4.12}

Schultz recently generalized in \cite{Schu} these results by
allowing $\Gal(K/k) \cong \mathds Z/p^n\mathds Z$ for any $n \in
\mathds N$, and removing the condition of cyclicity (as a module)
for the kernel. Shirbisheh also considers non cyclic kernels in
\cite{Shi}, where he studies embedding problems over the field
$\mathds Q(\zeta_{p^2})/\mathds Q(\zeta_p)$. Schultz's approach
differs in that he gives explicit descriptions for all possible
extensions of $\mathds Z/p^n\mathds Z$ by a finite $\mathds
F_p[G]$-module $A$, and then he finds a parameterizing set for each
such group within $J(K) = K^\times/K^{\times p}$. Note that the
study of the $\mathds F_p[\Gal(K/k)]$-structure of $J(K)$ was
initiated by Waterhouse in \cite{Wa}, and Schultz's results from
\cite{Schu} can be thought of as a completion of the ideas that
Waterhouse presents there.

We are going now to define three symbols that appear in the
statement of the main result from \cite{Schu}. For a field extension
$K/k$ with $\Gal(K/k) \cong \mathds Z/p^n\mathds Z$, let $K_i$
denote the intermediate field of degree $p^i$ over $k$.  If the
embedding problem $(K/k,\mathds Z/p^{n+1}\mathds Z,\mu_p)$ has a
solution, then define $i(K/k) = -\infty$.  Otherwise, let $s$ be the
minimum value such that the embedding problem $(K/K_s,\mathds
Z/p^{n-s+1}\mathds Z,\mu_p)$ has a solution, and define $i(K/k) =
s-1$.  Notice that we have $i(K/k) \in \{-\infty,0,\cdots,n-1\}$
provided $K \ne k$.

Denote by $\binom{n}{m}_p$ the $p$-binomial coefficient, defined for
$n \in \mathds N$ and satisfying
$$
\binom{n}{m}_p = \left\{\begin{array}{ll}
0, &\mbox{ if }m<0 \mbox{ or }m>n\\
\frac{(p^n-1)\cdots(p^{n-m+1}-1)}{(p^m-1)\cdots (p-1)},&\mbox{ if
}0\leq m \leq n.\end{array}\right.
$$

Finally, denote by $\lceil x\rceil$ the standard ceiling function,
i.e., the smallest integer not less than $x$.

\newtheorem{t4.13}[l4.1]{Theorem}
\begin{t4.13}\label{t4.13}
{\rm(\cite[Theorem 1]{Schu})} Let $G = \langle \sigma \rangle \simeq
\mathds Z/p^n\mathds Z$, where $p$ is a prime and $n>1$ when $p=2$,
and suppose that $K/k$ is an extension of fields so that $\Gal(K/k)
\simeq G$ and $\zeta_p \in K$.  Suppose that $A \simeq
\bigoplus_{i=1}^{p^n}\oplus_{d_i} \mathds F_p[G]/(\sigma-1)^i$ as an
$\mathds F_p[G]$-module, and write $\Delta(A_{\{i\}}) = \sum_{j \geq
i} d_j$.  For $1 \leq i \leq p^n$ let $$\mathfrak{D}_{\{i\}} =
\dim_{\mathds F_p}\left(\frac{N_{K_{\lceil
\log_p(i)\rceil}/k}(K_{\lceil \log_p(i)\rceil}^\times)}{K^{\times
p}}\right).$$ Then the embedding problem $(K/k,A \rtimes G,A)$ has a
solution over $K/k$ if and only if $\Delta(A_{\{i\}}) \leq
\mathfrak{D}_{\{i\}}$.

If $k^\times/K^{\times p}$ is infinite and the embedding problem
$(K/k,A \rtimes G,A)$ is solvable, then there are infinitely many
solutions to this embedding problem over $K/k$.  If
$k^\times/K^{\times p}$ is finite and the embedding problem $(K/k,A
\rtimes G,A)$ is solvable, then the number of solutions to this
embedding problem over $K/k$ is
$$\prod_{i=1}^{p^n}\binom{\mathfrak{D}_{\{i\}} - \Delta(A_{\{i+1\}}) - \mathds{1}_{i=p^{i(K/k)}+1}}{\Delta(A_{\{i\}})-\Delta(A_{\{i+1\}})}_p~p^{d_i \left(\sum_{j<i} \mathfrak{D}_{\{j\}} - \Delta(A_{\{j\}}) - \mathds{1}_{j=p^{i(K/k)}+1}\cdot \mathds{1}_{i=p^n})\right)}.$$
\end{t4.13}

\medskip

{\bf Example.} Let us consider the group $G_7=(\mathds Z/p\mathds
Z)^3\rtimes \mathds Z/p\mathds Z$ of order $p^4$ (see the definition
of $G_7$ after Proposition \ref{p4.10}). We may put $A= \mathds
F_p[G]/(\sigma-1)^3$, a cyclic module of dimension $3$, and
$G=\mathds Z/p\mathds Z$. Thus $G_7\cong A\rtimes G$, and we can
apply Theorem \ref{t4.13}.

We have $d_3=1$ and $d_j=0$ for $j\neq 3$, so $\Delta(A_{\{i\}}) =
0$ if $i\geq 4$ and $\Delta(A_{\{i\}}) = 1$ if $i \leq 3$. Therefore
$(K/k,A \rtimes G,A)$ is solvable if and only if
$\mathfrak{D}_{\{3\}} \geq 1$. As $\lceil\log_p(3)\rceil=1$ for
$p>2$, one simply need to check that $N_{K_1/k}(K_1^\times)$ is not
contained in $K^{\times p}$.\qed

\medskip

Schultz also investigates in \cite{Schu} the case when $K$ does not
contain a primitive $p$-th root of unity. In \cite[Section 6]{Schu}
$\mathfrak{D}_{i}$ is redefined in this case, and \cite[Theorem
6.3]{Schu} is the generalization of \cite[Theorem 1]{Schu} which
drops the Kummer theory assumption. Thus the example given above can
also be extended in the general case.

\section{Automatic realizations}
\label{8}

Let $G$ and $H$ be finite groups. If any field $k$ admitting a
$G$-extension also admits an $H$-extension, we will write
$G\Longrightarrow H$. A statement $G\Longrightarrow H$ is called an
\emph{automatic realization}. For instance, it is well-known (see
\cite{Wh}) that $\mathds Z/4\mathds Z\Longrightarrow \mathds
Z/2^n\mathds Z$ for all $n\in\mathds N$. Of course, we always have
the trivial automatic realization $G\Longrightarrow G/N$. Also, it
is not hard to verify that a necessary condition for the automatic
realization $G\Longrightarrow H$ to hold is the minimum number of
generators for $H$ is less than or equal to the minimum number of
generators for $G$ (see \cite[p. 1268]{Le-95}).

We keep the notations from Section \ref{5} for the four non abelian
$2$-groups having a cyclic subgroup of index $2$. Jensen and Yui
determined in \cite[Theorem III.3.6]{JY}) the automatic realization
$Q_8\Longrightarrow D_8$. Ledet proved in \cite[Proposition 5.8,
Proposition 5.10]{Le-95} the automatic realizations
$Q_{16}\Longrightarrow D_{16}, SD_{16}\Longrightarrow M_{16}$ and
many more. He also proved in \cite{Le-98} the automatic realization
$Q_{32}\Longrightarrow D_{32}$. Other automatic realizations
concerning groups of order $16$ are found by Grundman and Smith in
\cite{GS-16}, where they consider fields with special properties
(e.g. fields with a given level or a number of square classes).
There are as well a number of non trivial automatic realizations
among groups of order $32$ considered in \cite{Mi-32,GS-32,Sm}. More
automatic realizations among $2$-groups can be found in three papers
written by Jensen \cite{Je1,Je2,Je3} and in a paper written by Gao,
Leep, Min\'ac and Smith \cite{GLMS}.

We proceed with results about $p$-groups for an odd prime $p$.
Whaples showed in \cite{Wh} the automatic realization $\mathds
Z/p\mathds Z\Longrightarrow \mathds Z/p^n\mathds Z$ holds for all
$n\in\mathds N$. Brattstr\"om verified in \cite[Theorem 2]{Br} the
automatic realization $G_1\Longrightarrow G_2$, where $G_1$ and
$G_2$ are the two-non abelian groups of order $p^3$ defined in
Section \ref{7} (above Theorem \ref{t4.1}). Michailov showed in
\cite[Theorem5.2]{Mi-p4} that the automatic realization
$G_3\Longrightarrow G_4$ holds, where $G_3$ and $G_3$ are non
abelian groups of order $p^4$ defined in Section \ref{7} (below
Theorem \ref{t4.2}). In \cite{Br,Mi-p4} is shown also that the
reverse automatic realisations $G_2\Longrightarrow G_1$ and
$G_4\Longrightarrow G_3$ are not valid.

Min\'ac, Schultz and Swallow established in \cite{MSS} automatic
realizations among the semidirect products $M\rtimes G$, where
$G=\langle\sigma\rangle$ is a cyclic group of order $p^n$, and $M$
is a quotient of the group ring $\mathds F_p[G]$. For the group ring
$\mathds F_p[G]$ there exist precisely $p^n$ nonzero ring quotients,
namely $M_j:=\mathds F_p[G]/\langle(\sigma-1)^j\rangle$ for
$j=1,2,\dots,p^n$. Multiplication in $\mathds F_p[G]$ induces an
$\mathds F_p[G]$-action on each $M_j$. In particular, each $M_j$ is
a $G$-module.

\newtheorem{t8.1}{Theorem}[section]
\begin{t8.1}\label{t8.1}
{\rm(\cite[Theorem 1]{MSS})} $M_{p^i+c}\rtimes G\Longrightarrow
M_{p^{i+1}}\rtimes G$ for $0\leq i< n, 1\leq c<p^{i+1}-p^i$.
\end{t8.1}

There is a more general approach to the automatic realizations. We
write $\nu(G,k)$ for the number of distinct $G$-extensions of $k$
within a fixed algebraic closure of $k$.  We then write
$\mathfrak{K}(G)$ for the set of fields $k$ such that $\nu(G,k) \geq
1$. Now it is clear that if $k\in \mathfrak{K}(G)$ implies $k\in
\mathfrak{K}(H)$, then the automatic realization $G\Longrightarrow
H$ holds.

Berg and Schultz recently considered in \cite{BS} a close cousin of
automatic realization results.  The realization multiplicity of $G$,
written $\nu(G)$, is defined as $$\nu(G) = \min_{k \in
\mathfrak{K}(G)}\{\nu(G,k)\}.$$

\newtheorem{t8.2}[t8.1]{Theorem}
\begin{t8.2}\label{t8.2}
{\rm(\cite[Theorem 1.1]{BS})} Suppose that $p$ is prime and $n$ is a
positive integer, with $n \geq 2$ when $p=2$.  Let $k$ be given.
Then
$$\nu\left(\mathds F_p[\mathds Z/p^n\mathds Z]^k\rtimes \mathds Z/p^n\mathds Z\right) \geq p^k.$$
\end{t8.2}

Schults generalized the latter result in \cite{Schu}.

\newtheorem{t8.3}[t8.1]{Theorem}
\begin{t8.3}\label{t8.3}
{\rm(\cite[Theorem 1.3]{Schu})} Let $G = \langle \sigma
\rangle\simeq \mathds Z/p^n\mathds Z$, with $n>1$ when $p=2$.
Suppose that $A$ is an $\mathds F_p[G]$-module which is not
isomorphic to
$$\mathds F_p[G]/(\sigma-1)^{p^j+1} \bigoplus_{i=0}^{p^n} \oplus_{d_i}
\mathds F_p[G]/(\sigma-1)^{p^i}$$ for any choice of $j \in
\{-\infty,0,\cdots,n-1\}$ and $d_i \in \mathds Z$.  If $\hat G$ is
any extension of $G$ by $A$, and if $A$ contains elements
$a_1,\cdots,a_k$ which are $\mathds F_p[G]$-independent, then
$\nu(\hat G) \geq p^k$.
\end{t8.3}

\section*{Acknowledgements}
\label{9}

We are grateful to Prof. A. Schultz who read our manuscript, offered
constructive comments and provided the example at the end of Section
\ref{7}.

\end{document}